\documentclass[11pt]{amsart}

\usepackage{amsmath,amsfonts,amsthm,amssymb,
amsthm,graphicx,mathtools,amscd, enumerate,paralist,bm}
\usepackage{hyperref}

\usepackage[mathscr]{eucal}
\usepackage{graphics, color, epsfig}
\usepackage[nocompress]{cite}

\theoremstyle{plain}
\newtheorem{theorem}{Theorem}[section]
\newtheorem{lemma}[theorem]{Lemma}
\newtheorem{corollary}[theorem]{Corollary}
\newtheorem{definition}[theorem]{Definition}
\newtheorem{example}[theorem]{Example}
\newtheorem{proposition}[theorem]{Proposition}

\newtheorem{assumption}[theorem]{Assumption}

\newtheorem{remark}[theorem]{Remark}

\newcommand{\la}{\lambda}
\newcommand{\eps}{\varepsilon}
\newcommand{\ph}{\varphi}

\newcommand{\al}{\alpha}

\newcommand{\gam}{\gamma}

\newcommand{\sig}{\sigma}
\newcommand{\del}{\delta}
\newcommand{\delh}{\hat{\delta}}
\newcommand{\om}{\omega}

\newcommand{\Del}{\mathnormal{\Delta}}

\newcommand{\Om}{\mathnormal{\Omega}}

\newcommand{\N}{{\mathbb N}}

\newcommand{\R}{{\mathbb R}}

\newcommand{\EE}{{\mathbb E}}
\newcommand{\E}{{\mathbb E}}

\newcommand{\PP}{{\mathbb P}}

\newcommand{\XX}{{\mathbb X}}
\newcommand{\YY}{{\mathbb Y}}

\newcommand{\ONE}{\boldsymbol{1}}
\newcommand{\half}{{\textstyle\frac{1}{2}}}

\newcommand{\calB}{{\mathcal B}}
\newcommand{\calC}{{\mathcal C}}
\newcommand{\calD}{{\mathcal D}}

\newcommand{\calF}{{\mathcal F}}

\newcommand{\calH}{{\mathcal H}}
\newcommand{\calI}{{\mathcal I}}

\newcommand{\calM}{{\mathcal M}}

\newcommand{\calS}{{\mathcal S}}

\newcommand{\calV}{{\mathcal V}}

\newcommand{\bF}{{\mathbf F}}

\newcommand{\bX}{{\mathbf X}}

\newcommand{\frp}{\mathfrak{p}}

\renewcommand{\proof}{\noindent{\bf Proof.\ }}

\newcommand{\lan}{\langle}
\newcommand{\ran}{\rangle}

\newcommand{\supp}{{\rm supp}}

\newcommand{\w}{\wedge}

\newcommand{\pl}{\partial}

\newcommand{\iy}{\infty}

\newcommand{\loc}{{\rm loc}}

\newcommand{\noi}{\noindent}

\newcommand{\LE}{\preccurlyeq}

\newcommand{\bi}{\begin{itemize}}
\newcommand{\ei}{\end{itemize}}

\newcommand{\bb}{\zeta}

\newcommand{\bfx}{\mathbf{x}}
\newcommand{\bfX}{\mathbf{X}}

\newcommand{\atlas}{\textsf{Atlas}}
\newcommand{\ppp}{\textsf{PPP}}
\newcommand{\poi}{\textsf{Poi}}
\newcommand{\Exp}{\textsf{Exp}}

\newcommand{\pol}{\text{\rm pol}}

\newcommand{\BL}{\text{\tiny BL}}
\newcommand{\Lip}{\text{\tiny Lip}}

\newcommand{\one}{\mathbf{1}}

\newcommand{\cc}[1]{[#1]}

\usepackage{tcolorbox}
\newtcolorbox{bx}{colback=white!5!white,colframe=blue!75!black}

\definecolor{stred}{rgb}{1.0, 0.2, 0.37}

\begin{document}

\title{On free boundary problems for the Atlas model
}

\date{\today}

\author{Rami Atar}

\author{Amarjit Budhiraja}

\address{\small Viterbi Faculty of Electrical and Computer Engineering,
Technion -- Israel Institute of Technology, Haifa, Israel
\\
{\tt rami@technion.ac.il}
}

\address{\small Department of Statistics and Operations Research,
University of North Carolina at Chapel Hill, USA
\\
{\tt budhiraj@email.unc.edu}
}

\begin{abstract}
For $n\in\N$, let $\{X^n_i\}$ be an infinite collection of Brownian  particles on the real line where the leftmost particle $\min_iX^n_i(t)$ is given a drift $n$, and let $\mu^n_t=n^{-1}\sum_i\del_{X^n_i(t)}$, $t\ge0$ denote the normalized configuration measure.
The case where the initial particle positions follow a Poisson point process on $[0,\iy)$ of intensity $n\la$, $\la>0$ was studied in \cite{dem19}, where it was shown that $\mu^n_t$ converge, as $n\to\iy$, to a limit characterized by a Stefan problem of melting solid (respectively, freezing supercooled liquid) type when $\la\ge 2$ (respectively, $0<\la<2$). In this paper it is assumed that $\mu^n_0\to\mu_0$ in probability, where $\mu_0$ is supported on $[0,\iy)$ and satisfies a polynomial growth condition. Because $(y-x)^{-1}\mu_0((x,y])$, $0<x<y$
need not be bounded below or above by $2$, the model does not give rise to a Stefan problem of either of the above types. Under mild assumptions, it is shown that $\mu^n_t$ converge to a limit characterized by a free boundary problem involving measures.
Under the additional assumption that $\mu_0(dx)\ge\la_0\,{\rm leb}_{[0,\iy)}(dx)$ for some $\la_0>0$, the free boundary exists as a continuous trajectory, and the process determined by the leftmost particle converges to it.
\end{abstract}

\maketitle

\section{Introduction}

\subsection{The Atlas model and the Stefan problem}
\label{sec:model}

The infinite Atlas model consists of an infinite collection of
particles on the real line driven by Brownian motions (BMs), where at any given
time, the particle that is leftmost is given a drift $a \in (0,\infty)$. The precise definition is formulated
as a weak solution to the system of stochastic differential equations (SDE) indexed by $\N_0:=\{0,1,2,\ldots\}$,
\begin{equation}\label{60}
dX_i(t) = a\ONE_{\{X_i(t) = \min_jX_j(t)\}} dt + dW_i(t), \quad X_i(0) = x_i\ge0,\quad i \in \N_0.
\end{equation}

\begin{definition}\label{def1}
Given a parameter $a\ge0$ and a probability law $\pi$ on $\R^\iy_+$ (where $\R_+=[0,\iy)$),
an {\rm $\atlas(a,\pi)$} is a pair
$(X=(X_i)_{i\in\N_0},W=(W_i)_{i\in\N_0})$ of continuous adapted processes
defined on a filtered probability space
$(\Om,\calF,\{\calF_t\},\PP)$, such that
\\
(i) $W_i$ are mutually independent standard
$\{\calF_t\}$-BMs,
\\
(ii) $\PP$-a.s., \eqref{60} holds; in particular,
the minimum in \eqref{60} is attained for all $t$,
\\
(iii) The law of $\bfx=(x_i)_{i\in\N_0}$ under $\PP$ is $\pi$.
\end{definition}

To discuss the well-posedness of the model, let
\[
\calI=\{({\bf c}=(c_i)_{i\in\N_0}\in\R^{\N_0}:
c_0\le c_1\le\cdots\}
\]
denote the set of nondecreasing sequences and let
$\calI_+=\{{\bf c}\in\calI:c_0\ge0\}$.
We will always assume
that $\pi$ is supported on $\calI_+$.
A sequence $(c_i)_{i\in\N_0}$ is said to be {\it rankable} if 
there exists a bijection $\tau:\N_0\to\N_0$ such that
$(c_{\tau(i)})\in\calI$. For the drift coefficient in \eqref{60}
to be well-defined, it suffices that $X=(X_i(t))_{i\in\N_0}$,
is rankable for all $t$ a.s., via a measurable bijection $\tau=\tau_t$.
It is known that if $\bfx = (x_i)_{i\in \N_0} \in \calH$ a.s., where
$$
\calH \doteq \Big\{{\bf c} \in \R^{\infty}: \sum_{i=0}^{\infty} e^{-\al (c_i\vee0)^2} < \infty \mbox{ for all } \al>0\Big\},
$$
then \eqref{60} possesses a unique-in-law weak solution, which,
in particular, is rankable for all $t$ \cite{AS, ichiba2013strong}.
Thus, when $\pi$ is supported on $\calI_+\cap\calH$,
the $\atlas(a,\pi)$ exists and is unique in law.
As a result, there also exists a process
$Y=(Y_i(t))_{i\in\N_0}$ such that $Y_0(t)\le Y_1(t)\le\cdots$ and
$Y_i(t)=X_{\tau_t(i)}(t)$ for all $i$ and $t$ a.s.
The systems $X$ and $Y$ are referred to as the
{\it named} and, respectively, {\it ranked} particle system.

For $n\in\N$, let $(X^n,W)$ be a sequence
of $\atlas(a^n,\pi^n)$ models, where $a^n>0$,
$\pi^n$ is supported on $\calI_+\cap\calH$,
and throughout, the 
dependence of $W$ on $n$ is suppressed. Denote by
$\bfx^n$ the initial condition corresponding to $\pi^n$.
Without loss of generality, assume that
the entire sequence is defined
on a common filtered probability space
$(\Om,\calF,\{\calF_t\},\PP)$.
By the aforementioned result on weak uniqueness, the law of $X^n$,
under $\PP$ is uniquely determined by $(a^n,\pi^n)$
for each $n$. Following \cite{dem19}, we are interested
in a regime where the drift $a^n$ and the particle density both
grow with $n$. In particular, we take $a^n=n$, and,
denoting the normalized initial configuration measure by
$\mu^n_0=n^{-1}\sum_{i\ge0}\del_{x^n_i}$, where $(x^n_i)_{i \in \N_0} = \bfx^n$, assume
that $\mu^n_0$ converge to a locally finite measure $\mu_0$ on $\R_+:=[0,\iy)$, in the sense defined below.
Thus
\begin{equation}\label{61}
dX^n_i(t)=n\ONE_{\{X^n_i(t)=Y^n_0(t)\}}dt+dW_i(t),
\qquad
X^n_i(0)=x^n_i,\qquad i\in\N_0.
\end{equation}

Denote
\[
\mu^n_t(dx)=\frac{1}{n}\sum_{i\ge0}\del_{X^n_i(t)}(dx),
\qquad
\beta^n(dx,dt)=\del_{Y^n_0(t)}(dx)dt,
\qquad (x,t)\in\R\times\R_+.
\]
To define a suitable space for the measure-valued processes $\mu^n$,
let $\calM_*=\calM_*(\R)$
be the space of nonnegative measures $\mu$ on $(\R, \calB(\R))$ satisfying
$$\mu(-\infty, x]<\infty \mbox{ for all } x \in \R,$$
where, throughout, for a measure $\mu$ on $\R$,
$\mu(a,b]$ is shorthand notation for $\mu((a,b])$,
and a similar convention applies for intervals $(a,b)$, etc.
For $\mu\in\calM_*$, denote
\[
\bF^\mu(x)=\mu(-\iy,x], \qquad x\in\R.
\]
Equip the class $\calM_*$ with the following metric. 
Denote by $\text{BL}$ the class of real bounded Lipschitz functions $f$ on $\R$ such that $\text{supp}(f) \subset (-\infty, r]$ for some $r \in \R$. On this class
define the norm $\|\cdot\|_{\BL}$ as
$\|f\|_{\BL} \doteq \|f\|_{\infty} + \|f\|_{\Lip}$, $f \in \text{BL},$
where
$$\|f\|_{\infty} \doteq \sup_{x\in \R} |f(x)|, \qquad \|f\|_{\Lip} \doteq \sup_{x,y\in \R, x\neq y}\frac{|f(x)-f(y)|}{|x-y|}.$$
Equip $\calM_*$ with the metric
$$d_*(\mu, \nu) \doteq \sum_{r\in \N} 2^{-r} (1 \wedge d_{*,r}(\mu,\nu)),
\qquad \mu, \nu \in  \calM_*,
$$
defined via the pseudo-metric
\[
d_{*,r}(\mu,\nu)=\sup_{f\in \text{BL},\, \|f\|_{\BL} \le 1,\, \text{supp}(f) \subset (-\infty, r]} |\mu(f) - \nu(f)|, \qquad \mu, \nu \in  \calM_*,
\]
where, for a measure $\gamma$, $\gamma(f) = \int f d\gamma$.
Then $\calM_*$ is a Polish space, and $d_*$ is related to both
the weak and vague topologies; for example, if for some $r$,
$\nu^n\in\calM_*$ are measures supported on
$(-\iy,r]$ (respectively, $[r,\iy)$), then convergence of this sequence
in $(\calM_*,d_*)$ is equivalent to weak (respectively, vague) convergence.
Equip $C(\R_+,\calM_*)$ with the topology of uniform convergence
on compacts.
For each $n$, $\mu^n_t$ takes values in $C(\R_+,\calM_*)$,
as follows from \cite[Lemma 3.4]{AS} which shows that for all $t$,
$\bfX^n(t) = (X^n_i(t))_{i\in \N_0}$ takes values in $\calH$ and  is {\em locally finite} in the sense that a.s., for all $r \in \R$ and $T>0$, there exist finitely many $k$ such that $\min_{[0,T]} X^n_k \le r$.

One way to obtain a sequence $\atlas(n,\pi^n)$
is to consider the
$\atlas(1,\pi)$ model, $(X,W)$, with initial condition
$\bfx$, and rescale it as $X^n_i(t)=n^{-1}X_i(n^2t)$.
Because $n^{-1}W_i(n^2t)$ are standard BMs,
the rescaled processes form a sequence $\atlas(n,\pi^n)$,
where, under $\pi^n$, $(x^n_i)$ are equal in law to $(n^{-1}x_i)$.
This setting was studied in \cite{dem19} for
$\pi$ a Poisson point process (PPP) on $\R_+$ of intensity $\la\in(0,\iy)$ (i.e., $\pi^n$ is a PPP of intensity $n\la$).
It was shown that $\{\mu^n_t\}$ converge in probability
in $C(\R_+,\calM_*)$ to an $\calM_*$-valued path $\{\mu_t\}$, with $\mu_t$ supported on $[\sig_t,\iy)$, having density
$u(\cdot,t)$, where $(u,\sig)$ is the unique
solution to the following free boundary problem (FBP), which is an instance of the Stefan problem:
\begin{equation}\label{STE}\tag{STE}
\begin{split}
&(\pl_t-\half\pl_{xx})u=0,\hspace{3.9em} x>\sig_t,\ t>0,
\\
&u(\sig_t+,t)=2, \hspace{5.6em} t>0,
\\
&\textstyle\frac{d}{dt}\sig_t=-\frac{1}{4}\pl_xu(\sig_t+,t), \hspace{1.9em} t>0
\\
& \lim_{t\to0}u(x,t)=\la 1_{\{x>0\}}, \hspace{2.05em} x>0,
\end{split}
\end{equation}
extended by setting $u=0$ for $x\le\sig_t$. Here the solution is understood in a weak sense (see \cite[Definition 3.10]{dem19}). An explicit formula for the solution was found, and moreover, the large time asymptotics of $Y_0$ of the unscaled model was identified in terms of $\sig_t$ (see \cite[Theorem 1.2 and Corollary 1.4]{dem19}). In case where $\la\ge 2$, the free boundary $\sig_t$ is decreasing, whereas in the case $0<\la<2$ it is increasing. In the literature on the Stefan problem these cases are referred to as the melting of solid and, respectively, the freezing of supercooled liquid (abbreviated below to {\it melting} and {\it supercooled} respectively). Whereas the problem is well-posed in the melting case, the supercooled case is ill-posed: In some instances of this problem, classical global solutions do not exist, as singularities may appear in finite time; moreover, the proof of uniqueness may differ depending on the details of the problem. For this reason, the melting and supercooled cases were treated separately in \cite{dem19} as far as uniqueness of weak solutions is concerned, and moreover, in the latter case, uniqueness was proved in a restricted function class, tailored to the particular details of the problem. For another example of a particle system characterized at the hydrodynamic limit by a Stefan problem, with uniqueness in the supercooled case requiring a restricted function class, see \cite{chayes1996hydrodynamic}.

Our goal here is to extend the convergence result of \cite{dem19} to initial configurations that need not form a PPP and, more importantly, for which the limiting measure $\mu_0$ need not have a constant density (in fact, it may not be absolutely continuous w.r.t.\ the Lebesgue measure). When $\mu_0$ has a density bounded below or above by $2$, one still expects the limit to behave according to the melting or, respectively, the supercooled Stefan problem. However, in general, the solution will not possess a monotone free boundary, and a different set of tools will be required.

Our approach is based on a measure-theoretic notion of a FBP introduced in \cite{ata25}, where it was used to characterize the hydrodynamic limit of a branching Brownian particle system with selection. An advantage of this formulation is that it enables a treatment that does not require free boundary monotonicity, and in fact does not necessitate the existence of a free boundary as a trajectory.

The study of relations between particle system dynamics and parabolic FBP has been an active research field in the last decade. In \cite{dembo2019criticality}, the asymptotics of a particle system which can be viewed as a discrete version of the Stefan problem were addressed. In \cite{de-masi-excl}, a simple exclusion process with boundary behavior was studied at the hydrodynamic limit. Brownian particle systems with selection giving rise to FBP at the hydrodynamic limit were studied in one dimension \cite{de-masi-book, de-masi-nbbm, ber19, ata25}, and in general dimension \cite{ber21, ber22bee}. Integro-differential FBP were obtained for particle systems with non-local branching \cite{dur11, de-masi-nbbm2, ata-dr}.

In recent years, relaxed solutions to the supercooled Stefan problem have been introduced in \cite{delarue22}, via a probabilistic formulation, and in \cite{kim24, choi24, ghouss19}, based on a formulation involving measures. The measure-theoretic formulation of \cite{kim24, choi24, ghouss19} is different than the one we use in this paper (for example, in our version, there is no analogue to the $\nu$ appearing in \cite{choi24}, and in \cite{choi24} there is no analogue of our requirement that $\beta_t$, $t\in\R$ are probability measures; see below). Moreover, their main interest is in existence of solutions for the problem in general dimension (though their results are new even in dimension 1), and they address the melting and supercooled cases, whereas in this paper the treatment of non-monotone boundary is of crucial importance.

Relaxed solutions to a certain FBP have also been proposed in the aforementioned \cite{de-masi-book}, in the context of a particle systems with selection.

To define the FBP notion used in this paper, and state the main results, we need to introduce some additional notation.

\subsection{Notation}

Let $\calM^{(1)}(\R\times\R_+)$ denote the space of nonnegative Borel measures $\beta$ on $\R\times\R_+$ satisfying $\beta(\R\times[0,t])=t$ for all $t$, endowed with the topology of weak convergence on $\R\times[0,t]$ for every $t$. For $\beta$ in this class denote its disintegration as $\beta(dx,dt)=\beta_t(dx)dt$, and note that $\beta_t$ is a probability measure on $\R$ for a.e.\ $t$.

Denote by $\calM_*^{\pol}$ the set of $\mu\in\calM_*$
satisfying $\mu(-\iy,x]\le c(1+x_+^c)$ for some $c>0$,
and $\mu(\R)=\iy$.
Denote by $\calV^\pol(\R)$ the set of nondecreasing
right-continuous functions $v:\R\to \R_+$
satisfying $\lim_{x\to-\iy}v(x)=0$, $\lim_{x\to\iy}v(x)=\iy$ and
$v(x)\le c(1+x_+^c)$ for some $c$.
Then $v=\bF^\mu(\cdot)=\mu(-\iy,\cdot]$ is a bijection
from $\calM_*^\pol$ to $\calV^\pol(\R)$.
Let $\calV^\pol_0(\R)$ (respectively, $\calV^\pol_\#(\R)$) be the set of $v\in\calV^\pol(\R)$
for which $v=0$ on $(-\iy,0)$ (respectively,
$v|_{(-\iy,0)}\in L^1(-\iy,0)$).

Next, let $C^\pol(\R)=\calV^\pol(\R)\cap C(\R)$ and
$C^\pol_\#(\R)=\calV^\pol_\#(\R)\cap C(\R)$.
Finally, let $C^\pol_\loc(\R\times(0,\iy))$
(respectively, $C^\pol_{\#,\loc}(\R\times(0,\iy))$)
denote the set of functions $v\in C(\R\times(0,\iy))$
such that for every $t>0$,
$v(\cdot,t)\in C^\pol(\R)$ (respectively, $C^\pol_\#(\R)$), and
there exists $c=c(t)\in(0,\iy)$ such that
\[
\sup_{s\in(0,t]}v(x,s)\le c(1+x_+^c), \qquad x\in\R.
\]

For functions $f,g: \R \to \R$, we write
$\lan f,g\ran=\int fgdx$, and for a real function $f$ and a measure $\beta$ on a measurable space, we write $\lan f,\beta\ran=\int fd\beta$.

\subsection{FBP involving measures}

Whereas the function $u$ in equation \eqref{STE} represents density, we will be concerned with an equation for the cumulative mass function $v$, related to $u$ via $v(x,t)=\int_{-\iy}^xu(y,t)dy$.
Given $v_0 \in \calV_0^\pol(\R)$, representing initial cumulative distribution of mass, consider
\begin{equation}\label{FBP}\tag{FBP}
\begin{split}
&(\pl_t-\half\pl_{xx})v=0,\qquad x>\sig_t,\ t>0,
\\
&v(\sig_t+,t)=0,\qquad \pl_xv(\sig_t+,t)=2,\qquad t>0,
\\
& \lim_{t\to0}v(x,t)=v_0(x) \qquad \text{at all continuity points $x\ge 0$
of $v_0$}.
\end{split}
\end{equation}
Note that formally differentiating $\frac{d}{dt}v(\sig_t,t)=0$ in \eqref{FBP} gives $\frac{d}{dt}\sig_t=-\frac{1}{4}\pl_{xx}v(\sig_t,t)$, in agreement with \eqref{STE}.

\begin{definition}[Weak solution to \eqref{FBP}]\label{def00}
Given $v_0 \in \calV_0^\pol(\R)$, a pair $(v,\sig)\in C^\pol_\loc(\R\times(0,\iy))\times C((0,\iy),\R)$
is a weak solution to \eqref{FBP} if
for $\ph\in C^\iy_c(\R\times\R_+)$,
\begin{equation}\label{03-}
-\int_0^\iy\lan\pl_t\ph+\half\pl_{xx}\ph,v\ran dt
=\lan\ph(\cdot,0),v_0\ran-\int_0^\iy\ph(\sig_t,t)dt,
\end{equation}
and $v$ vanishes on $\{(x,t):x\le\sig_t,t>0\}$.
\end{definition}
For a measure-theoretic formulation of \eqref{FBP}, consider
\begin{equation}\label{MFBP}\tag{MFBP}
\begin{split}
&{\it i.} \hspace{4em} (\pl_t-\half\pl_{xx})v=-\beta,
\\
&{\it ii.} \hspace{3.9em} v(\cdot,0)=v_0,
\\
&{\it iii.} \hspace{3.1em} \int_{\R\times(0,\iy)}v\,d\beta=0.
\end{split}
\end{equation}

\begin{definition}[Solution to \eqref{MFBP}]\label{def0}
Given $v_0 \in \calV_0^\pol(\R)$, a pair $(v,\beta)\in C^\pol_\loc(\R\times(0,\iy))\times\calM^{(1)}(\R\times\R_+)$
is a solution to \eqref{MFBP} if
for $\ph\in C^\iy_c(\R\times\R_+)$,
\begin{equation}\label{03}
-\int_0^\iy\lan\pl_t\ph+\half\pl_{xx}\ph,v\ran dt
=\lan\ph(\cdot,0),v_0\ran-\int_{\R\times\R_+}\ph\, d\beta,
\end{equation}
and (\ref{MFBP}.iii) holds.
\end{definition}

If $(v,\sigma)$ is a weak solution of \eqref{FBP} then, with
\begin{equation}\label{w1}
\beta(dx,dt)=\del_{\sigma_t}(dx)dt,
\end{equation}
$(v,\beta)$ is a solution to \eqref{MFBP} (equation \eqref{03} follows from \eqref{03-}, and (\ref{MFBP}.iii) follows from the fact that $v$ vanishes on the curve $x=\sig_t$).

We will not be concerned with classical solutions of \eqref{FBP}, only with weak solutions of \eqref{FBP} and solutions to \eqref{MFBP}. However, for completeness, we mention in the appendix
well known conditions under which a classical solution $(v,\sig)$ of \eqref{FBP} is a weak solution, and consequently gives rise to a solution of \eqref{MFBP} via the same substitution \eqref{w1}.

\subsection{Main results}
Throughout we consider $\atlas(n,\pi^n)$, and $\mu^n, \beta^n$ as defined in Section \ref{sec:model}.
\begin{assumption}\label{ass1}
(a)
There exists $\mu_0\in\calM_*^\pol$ supported in $\R_+$,
such that $\mu^n_0\to\mu_0$ in probability in $(\calM_*,d_*)$.
\\
(b)
There exist constants $m\ge0$, $C, c>0$ such that
\begin{equation}\label{70}
\PP\left( \frac{\mu^n_0[j, j+1]}{1+j^m} >y\right) \le Ce^{-cy}, \mbox{ for all } n \in \N,\,j \in \N_0,\, y>0.
\end{equation}
\end{assumption}

\begin{remark}
Under Assumption \ref{ass1},
$\bfx^n\in\calH$ a.s.\ for each $n$. To see this,
fix $n$ and note, by Borel-Cantelli, that there exists
a (random) polynomial $p(\cdot)$ such that $\mu^n[0,j]\le p(j)$ for all $j$, a.s., and as a result, there exists a (random) $c>0$ such that
$x^n_i\ge ci^c$ for all $i$ a.s.
\end{remark}

Following are two important cases where Assumption \ref{ass1} holds.

\begin{example}[PPP initial condition]\label{ex1}
Suppose that, for each $n\in\N$, $\bfx^n$ forms a Poisson point process
on $\R_+$ with locally finite intensity measure $\gamma^n$, and that
$n^{-1}\gamma^n\to\mu_0$ in $(\calM_*,d_*)$,
where $\mu_0\in\calM^\pol_*$.
Suppose, moreover, that $\theta\in\calV^\pol_0(\R)$, where
$\theta(x)=\sup_n n^{-1}\gamma^n[0,x]$.
Then Assumption \ref{ass1} holds.
This is proved below in Proposition \ref{prop01}. Note that the setting considered in \cite{dem19} corresponds to $\gamma^n(dx) = n\lambda dx$ which clearly satisfies this sufficient condition. More generally, if $\gamma^n(dx) = n \la(x) dx$ where $\la(\cdot)$ grows at most polynomially and $\int_0^{\infty} \la(x) dx = \infty$, then $\gamma^n$ satisfies this condition. As another example, let $\la: \R_+ \to \R_+$ be a bounded function such that $\la(x) \to \la_0$ as $x \to \infty$. Then $\gamma^n(dx) = \la^n(x) dx$, where $\la^n(x) = n \la(nx)$ satisfies the assumption with $\mu_0(dx) = \la_0 dx$.
\end{example}

\begin{example}[Deterministic initial condition]\label{ex2}
Suppose $x^n_i=f(i/n)$ for some strictly increasing $f:\R_+\to\R_+$, $f(0)=0$, with inverse $f^{-1}\in\calV_0^\pol(\R)$.
Then Assumption \ref{ass1} holds with
$\mu_0[0,x]=f^{-1}(x)$. Part (a) of the assumption  follows on noting that, for all $j \in \N_0$,
$$|\mu^n_0[0,j] - \mu_0[0,j]| = |\mu^n_0[0,j] - f^{-1}(j)| \le 1/n,$$
whereas part (b) is immediate from the definition of $\{x^n_i\}$.
\end{example}

Our first main result shows that the limit of $(\mu^n,\beta^n)$ exists and can be characterized via \eqref{MFBP}.

\begin{theorem}\label{th1}
(a) Let $v_0\in\calV_0^\pol(\R)$. Then there exists a unique solution $(v,\beta)$ in the space
$C^\pol_{\#,\loc}(\R\times(0,\iy))\times \calM^{(1)}(\R\times\R_+)$ to \eqref{MFBP}.
\\
(b) Let Assumption \ref{ass1} hold and denote $v_0=\bF^{\mu_0}$.
Denote by $(v,\beta)$ the corresponding unique solution to \eqref{MFBP},
and, for $t>0$, let $\mu_t$ be defined via the relation
$v(\cdot,t)=\bF^{\mu_t}$.
Then one has $\{\mu_t\}_{t\ge0}\in C(\R_+,\calM_*)$. Moreover,
$(\mu^n,\beta^n)\to(\mu,\beta)$
in $C(\R_+,\calM_*)\times\calM^{(1)}(\R\times\R_+)$ in probability.
\end{theorem}

Given $(v,\beta)$ as above, one may regard the path
\begin{equation}\label{20d}
\sig_t=\inf\{x:v(x,t)>0\},\qquad t>0
\end{equation}
as a free boundary. However, this does not seem to be very useful in general because a priori there is no guarantee that the graph of $\sig$ equals the topological boundary of $\{(x,t):v(x,t)=0\}$ (equivalently, of $\{(x,t):v(x,t)>0\}$).
Our second main result gives conditions under which \eqref{20d} does determine the topological boundary, and relates the particle system's limit to \eqref{FBP}.
\begin{theorem}\label{th2}
Let $v_0\in\calV_0^\pol(\R)$ be such that, for some $\la_0>0$, $v_0(x) \ge \la_0 x$ for all $x\ge 0$. Let $(v,\beta)$ be as in Theorem \ref{th1}(a). Then the following hold.
\\
(a) $\beta$ takes the form $\beta(dx,dt)=\del_{\sig_t}(dx)dt$, where $\sig_t$ is as in \eqref{20d}. Moreover, defining $\sig_0=0$, one has $\sig\in C([0,\iy),\R)$.
\\
(b) The pair $(v,\sig)$ thus defined is the unique weak solution to \eqref{FBP}.
\\
(c) Suppose now that Assumption \ref{ass1} holds and, for each $n$, the initial gaps $\{Z^n_i(0)\}_{i \in \N}$, $Z^n_i(0) = x^n_i- x^n_{i-1}$, $i \in \N$, are stochastically dominated by $\{n^{-1}\bar Z_i\}_{i\in \N}$ where $\bar Z_i$ are i.i.d.\ exponentials of parameter $c$, for some $c>0$. Then $Y^n_0\to\sig$ in $C([0,\iy),\R)$, in probability.
\end{theorem}

\subsection{About the assumptions}

Although we have chosen to work under a polynomial growth condition, the results of this paper are mostly valid under exponential growth. In particular, the uniqueness results (Theorems \ref{th1}(a) and \ref{th2}(a,b)) will hold without any change to the proofs, whereas the for the convergence results (Theorems \ref{th1}(c) and \ref{th2}(c)), only minor modification of the assumptions and proofs are required to replace polynomial growth by an exponential growth.

We now go back to Example \ref{ex1}.

\begin{proposition}\label{prop01}
Under the setting of Example \ref{ex1},
Assumption \ref{ass1} holds.
\end{proposition}
\proof
To prove part (a) of the assumption, it suffices to show that, for
$f\in C_b(\R)$ supported on $(-\iy,r]$, for some $r$,
$B_n\to B$ in probability, where
$B_n=\lan f,\mu^n_0\ran$, $B=\lan f,\mu_0\ran$.
Now, $\E B_n=C_n:=n^{-1}\int f\gamma^n(dx)\to B$ by assumption.
Hence it suffices to prove that $\E[(B_n-C_n)^2]\to0$.
To this end, note that
\[
M^n_t:=\int_{[0,t]}f(x)\{\mu^n_0(dx)-n^{-1}\gamma^n(dx)\},
\qquad t\ge 0
\]
is a martingale, and
\[
[M^n]_r=n^{-2}\sum_{i:x^n_i\le r} f(x^n_i)^2
\le n^{-2}|f|^2_\iy N^n,
\]
where $N^n=\#\{i:x^n_i\le r\}\sim\text{\rm Poisson}(l^n)$,
$l^n=\gamma^n[0,r]$.
Hence
\[
\E[(B_n-C_n)^2]=\E[(M^n_r)^2]
\le c n^{-2}l^n.
\]
By assumption, $n^{-1}l^n$ is bounded.
This shows that $\E[(B_n-C_n)^2]\to0$ as required.

As for part (b) of Assumption \ref{ass1}, note
by the assumptions made in Example \ref{ex1},
that there exist $\al,m\in(0,\iy)$ such that
for all $n$ and $j$, $s_{n,j}:=\int_0^{j+1}\la^n(z)dz\le n\al(1+j^m)$. Fix such $\al,m$, and write $N(s)$ for a Poisson r.v.\ of parameter $s$.
Then, with $a=e^\al-1$,
the left-hand side of \eqref{70} is bounded, for all
$y>a$, by
\begin{align*}
\PP( N(s_{n,j})>(1+j^m)ny)
&\le e^{-\al(1+j^m)ny}e^{as_{n,j}}
\le e^{-\al(1+j^m) ny+an\al(1+j^m)}
\\
&\le e^{\al(-y+a)}=Ce^{-cy},
\end{align*}
where $C=e^{\al a}$, $c=\al$.
For $0\le y\le a$, the above bound is trivial.
This verifies part (b) of Assumption \ref{ass1}.
\qed

\begin{remark}
(a)
As already mentioned,
Assumption \ref{ass1} allows for settings where $\atlas(n,\pi^n)$ arises as the scaling of $\atlas(1,\pi)$ for some fixed $\pi$ supported on $\calI_+\cap\calH$.
One such example was discussed in Example \ref{ex1} where $\pi$ was the law of PPP on $\R_+$ with an intensity function $\la(\cdot)$ and $\la(x) \to \la_0$ as $x\to \infty$. More generally
consider the probability measure $\pi$ on $\calI_+$ given as the law of $(x_i)_{i \in \N_0}$ where for some $c, C, c_1>0$
\[
\PP\left(\left|\frac{\#\{i:x_i\le r\}}{r}- c_1\right| \ge y\right) \le C(e^{-cy} \wedge \al(r,y)) \text{ for all } r \ge r_0, \, y \ge 0,
\]
where $\al(r,y) \to 0$ as $r \to \infty$ for every $y>0$.
Then Assumption \ref{ass1} is satisfied. Indeed,
\[
\mu^n_0(-\iy,x]=n^{-1}\#\{i:\hat x^n_i\le x\}=n^{-1}\#\{i:x_i\le nx\}\to c_1x_+,
\]
in probability, and the first part of the assumption holds with the measure $\mu_0$ having
density $c_11_{[0,\iy)}(x)$. The second part of the assumption is easily verified with $m=1$ from the tail bound given above.
\noi

(b)  Example \ref{ex1} discussed settings with initial configurations given as PPP where the density of $\mu_0$ has polynomial growth. More generally, suppose that $\pi^n$  supported on $\calI_+$ is given as the probability law of the configuration $\{x^n_i\}$ such that, for some $p>0$,
for every $x\ge 0$,
\[
\frac{\#\{i:x^n_i\le x\}}{n}\to cx^p \text{ in probability, as $n\to\iy$,}
\]
and \eqref{70} holds. Then
\[
\mu^n_0(-\iy,x]=n^{-1}\#\{i:x^n_i\le x\}\to cx_+^p,
\]
and  Assumption \ref{ass1} holds  with $\mu_0$ having density
$cx^p1_{[0,\iy)}(x)$. An elementary example where these properties hold is $x^n_i=(i/cn)^{1/p}$, $i \in \N_0$.

\end{remark}

\subsection{Organization of the paper}
In Section \ref{sec2} we prove tightness of the sequence $(\mu^n,\beta^n)$. Section \ref{sec3} shows that limits are supported on solutions to \eqref{MFBP}, and
Section \ref{sec4} proves uniqueness of solutions to \eqref{MFBP}.
The results of these sections are then combined in Section \ref{sec5} to obtain the proof of the main results.

\section{Tightness}\label{sec2}
In this section we prove tightness of the pair $(\mu^n,\beta^n)$, needed in the proof of Theorem \ref{th1}(b).

We will use the following notation.
Denote the Euclidean norm in $\R^k$ by $\|\cdot\|$.
For $\bX$ a set equipped with a pseudo-metric $d$, $\xi:\R_+\to\bX$, and $0\le \del\le T<\iy$, let
\begin{align*}
w^{(d)}_T(\xi,\del)&=\sup\{d(\xi(s),\xi(t)):s,t\in[0,T],|s-t|\le\del\}.
\end{align*}
For $\xi:\R_+\to\R^k$, $\|\xi\|^*_T=\sup\{\|\xi(t)\|:t\in[0,T]\}$.
Let $W$ be a standard BM.

Let $\nu^{n}$ be a sequence of $C(\R_+,\calM_*)$-valued random variables.
The conjunction of \eqref{50} and \eqref{51} below
is a well-known sufficient condition for the tightness of the sequence
(cf. \cite[Theorem 14.5]{kal1},\cite[Theorem 3.7.2]{ethkur}).
\begin{equation}\label{50}
\text{
$\forall$ $\eta>0$ and $t\ge0$, $\exists$ a compact $\Gamma\in\calM_*$ such that}
\liminf_{n\to\iy}\PP(\nu^n_t\in\Gamma)\ge1-\eta,
\end{equation}
\begin{equation}\label{51}
\text{
$\forall$ $T>0$ and $\eta>0$, $\exists$ $\del>0$ such that, with $w_T=w^{(d_*)}_T$,}
\limsup_{n\to\iy}\PP(w_T(\nu^n,\del)\ge\eta)\le \eta.
\end{equation}
Now, if $M_n$, $a_n$ and $\eta_n$ are $(0,\iy)$-valued sequences with $\eta_n\to0$, and $a_n \to \infty$, the set
\[
\Gamma_{\{M_n,a_n, \eta_n\}}
=\{\mu\in\calM_*:\mu(-\iy,n]\le M_n,\ \mu(-\iy,-a_n]\le\eta_n, \ n\in\N\}
\]
is compact.
Therefore to verify \eqref{50}, it suffices to prove that
\begin{equation}\label{eq:bb1}
\lim_{M\to \infty} \limsup_{n\to \infty} \PP(\nu^n_t(-\iy,x] \ge M) = 0 \quad \mbox{ for every } x \in \R \mbox{ and } t\ge 0,
\end{equation}
\begin{equation}\label{eq:bb2}
\lim_{x\to -\infty} \limsup_{n\to \infty} \PP(\nu^n_t(-\iy,x] \ge\eta) = 0 \quad \mbox{ for every $t\ge 0$ and $\eta>0$}.
\end{equation}
To verify \eqref{51} it suffices that, for every $r\in\N$ and $\eta>0$ there exists $\del>0$
such that
\begin{equation}\label{52}
\limsup_n\PP(w^{(d_{*,r})}_T(\nu^n,\del)\ge\eta)\le\eta.
\end{equation}

\begin{lemma}\label{lem1}
Let Assumption \ref{ass1} hold. Then we have the following.
\begin{enumerate}[(a)]
\item For every $n$, $\mu^n$ has sample paths in $C([0,\infty),\calM_*)$.
\item Denote
\[
R^n(x,t)=\sum_{i\ge0}\PP(\min_{[0,t]}X^n_i\le x)^{1/2}, \qquad x\in\R,\, t\in\R_+.
\]
Then for every $x$ and $t$,
$\limsup_n n^{-1}R^n(x,t)<\iy$.
\item For every $\eta >0$ and $T>0$,
$$\lim_{x \to -\infty} \sup_n \PP(\sup_{t \le T} \mu^n_t(-\infty, x) > \eta) =0.$$
\item The sequence $\{\mu^n\}_{n\in \N}$ is tight
in $C([0,\infty),\calM_*)$.
\item If $\mu$ is a subsequential limit of $\{\mu^n\}$ then
for every $t$, $\int_{-\iy}^0\mu_t(-\iy,x]dx<\iy$,
and for every
$T$ there exists a (random) $c\in(0,\iy)$ such that
\[
\sup_{t\le T}\mu_t(-\iy,y]\le c(1+y^c),
\qquad y\ge0,
\]
\end{enumerate}
\end{lemma}

\proof
(a) Since, for each $n \in \N$, $\pi_n$ is supported on $\calH$, from  \cite[Lemma 3.4]{AS},
$\{X^n_k,k\in\N_0\}$ is {\em locally finite}, i.e. a.s., for all $r \in \R$ and $T>0$, there exist finitely many $k$ such that $\min_{[0,T]} X^n_k \le r$. Continuity of $\mu^n$ is immediate from this and the a.s. continuity of $X^n_k$ for each $n$ and $k$.

(b) Fix $x$ and $t$ and assume without loss of generality that
$x>0$. Fix $0<q<1$. Suppressing $x$ and $t$, we have $R^n \le R^n_1 + R^n_2$, where
	$$
	R^n_1 = \sum_{i\ge 0} \PP(x^n_i \le (i/n)^q)^{1/2}, \qquad
	R^n_2 = \sum_{i\ge 0} \PP\Big((i/n)^q + \min_{[0,t]} W(t) \le x\Big)^{1/2}.
	$$
To bound $R_1^n$, note that, for $a>0$, with $m$, $c$ and $C$ as in Assumption \ref{ass1},
	\begin{multline*}
		\PP(x^n_i\le a) = \PP\Big(\sum_{k} 1\{x^n_k \in [0,a]\} \ge i+1\Big)\\
		\le \sum_{j=0}^{\lfloor a\rfloor} \PP\left(\sum_{k} 1\{x^n_k \in [j,j+1]\} \ge \frac{i}{\lfloor a\rfloor +1}\right)\\
		= \sum_{j=0}^{\lfloor a\rfloor} \PP\left(\frac{\mu^n[j, j+1]}{1+j^m} \ge \frac{i}{n (\lfloor a\rfloor +1) (1+j^m)}\right)\\
		\le \sum_{j=0}^{\lfloor a\rfloor} \PP\left(\frac{\mu^n[j, j+1]}{1+j^m} \ge \frac{i}{n (\lfloor a\rfloor +1)^{m+1} }\right)\\
		\le (a+1)C \exp\left\{ - c\frac{i}{n (a +1)^{m+1}}\right\}.
		\end{multline*}
	Taking $a = (i/n)^q$, with $c_0 = c/2$,
	\begin{multline*}
	\frac{R^n_1}{n} = 	\frac{1}{n}\sum_{i\ge 0} \PP(x^n_i \le (i/n)^q)^{1/2}\\
	\le \frac{C}{n}\sum_{i\ge 0} ((i/n)^q+1)^{1/2} \exp\left\{ - c_0\frac{i}{n ((i/n)^q +1)^{m+1}}\right\}\\
	\le C\int_0^{\infty} (x^q+2)^{1/2} \exp\left\{-c_0 \frac{(x-1)^+}{(x^q+2)^{m+1}}\right\} dx.
	\end{multline*}
	Choose $q>0$ such that $q(m+1)<1$. Then we have that the last integral is finite, showing that
	$\limsup_{n\to \infty} n^{-1}R_1^n<\infty$.

To bound $R^n_2$, write
	$$R^n_2 \le (2x)^{1/q}n + \sum_{i\ge 0} \PP\Big(\|W\|^*_t \ge \frac{1}{2} (i/n)^q\Big)^{1/2}.$$
For some $c_1, c_2>0$
$$\PP\Big(\|W\|^*_t \ge \frac{i^q}{2n^q}\Big)^{1/2} \le c_1 e^{- c_2 i^{2q}/n^{2q}}.$$
Also
$$\sum_{i\ge0}e^{- c_2 i^{2q}/n^{2q}} \le (n+1)+\int_n^{\infty} e^{- c_2 y^{2q}/n^{2q}} dy.$$
Making the substitution $z= y^{2q}$, followed by the substitution
$u=c_2z/n^{2q}$, the above expression is bounded by
\begin{align}\label{68}
(n+1)+\frac{1}{2q}\int_{n^{2q}}^{\infty} z^{1/(2q)-1} e^{-c_2z/n^{2q}} dz
= (n+1)+\frac{1}{2q} nc_2^{-1/(2q)}\int_{c_2}^{\infty} u^{1/(2q)-1} e^{-u} du \le c_3 n,
\end{align}
where $c_3<\iy$.
This shows $\limsup_nn^{-1}R^n_2<\iy$.
The result follows.

(c) Fix  $\eta >0$ and $T> 0$.
It follows from Assumption \ref{ass1} that there exists
$c'$ such that, for all $n$ and $j$,
\begin{equation}\label{eq12.1}
\E\mu^n_0[j,j+1]\le c'(1+j^m).\end{equation} 
For $x\ge0$,
using independence of $\{x^n_i\}$ and $\{W_i\}$,
\begin{align}
\notag
\E\sup_{t \in [0,T]}\mu^n_t(-\iy,-x]&\le \frac{1}{n}\E\sum_{i,j}\sup_{t\in [0,T]}1_{\{X^n_i(t)\le-x,\,x^n_i\in[j,j+1)\}}\le\frac{1}{n}\sum_{i,j}\PP(\|W_i\|^*_t\ge x+j)\PP(x^n_i\in[j,j+1])
\\
&\le \sum_jc_1e^{-c_2(x+j)^2}\E[\mu^n_0[j,j+1]]
\le c'c_1\sum_je^{-c_2(x+j)^2}(1+j^m)
\le c_3e^{-c_2x^2},
\label{71}
\end{align}
for some constants $c_1, c_2, c_3<\iy$. Thus
$\limsup_n\PP(\sup_{t \in [0,T]}\mu^n_t(-\iy,-x]\ge\eta)\le \eta^{-1}c_3e^{-c_2x^2}$,
and sending $x\to\iy$ completes the proof of (c).

(d) We prove tightness of $\mu^n$, by verifying
\eqref{eq:bb1}, \eqref{eq:bb2} and \eqref{52}.
To show that $\mu^n$ satisfy \eqref{eq:bb1}, observe that
\[
\PP(\mu^n_t(-\iy,x]\ge M)\le\frac{1}{Mn}\EE\sum_{i\ge0}1_{\{X^n_i(t)\le x\}}
\le\frac{R^n(x,t)}{Mn}.
\]
By part (b), sending $n\to\iy$ and then $M\to\iy$
gives \eqref{eq:bb1}.

Next note that \eqref{eq:bb2}, is immediate from (c).

We next show that $\mu^n$ satisfy \eqref{52}.
Fix $r$ and $T$. Then, for $0\le s<t\le T$ and $f$ such that
$\supp(f)\subset(-\iy,r]$ and $\|f\|_\BL\le 1$,
\begin{align*}
|\mu^n_t(f)-\mu^n_s(f)|
&\le\frac{1}{n}\sum_{i\ge0}|X^n_i(t)-X^n_i(s)|1_{\{\min_{[s,t]}X^n_i<r\}}
\\
&\le \frac{1}{n}\sum_{i\ge0}\Big[|W_i(t)-W_i(s)|1_{\{\min_{[s,t]}X^n_i<r\}}
+n\int_s^t1_{\{X^n_i(u)=Y^n_0(u)\}}du\Big],
\end{align*}
where \eqref{61} is used. Hence
\[
w_T^{(d_{*,r})}(\mu^n,\del)
\le\frac{1}{n}\sum_{i\ge0}w_T(W_i,\del)1_{\{\min_{[0,T]}X^n_i<r\}}
+\del.
\]
Thus, denoting $\rho(\del)=\E[w_T(W,\del)^2]^{1/2}$,
\[
\E w_T^{(d_{*,r})}(\mu^n,\del)\le \frac{1}{n}\rho(\del)
\sum_{i\ge0}\PP(\min_{[0,T]}X^n_i<r)^{1/2}+\del
= \rho(\del)\frac{R^n(r,T)}{n}+\del.
\]
Using part (b), we have for some constant $c$,
\[
\limsup_n\E w_T^{(d_{*,r})}(\mu^n,\del)\le c\rho(\del)+\del.
\]
Above, the expression on the right converges to $0$ as $\del\to0^+$,
which, by Chebychev's inequality, shows that $\{\mu^n\}$ satisfies \eqref{52}.
This completes the proof of tightness of $\mu^n$.

(e)
For the first assertion, by \eqref{71} and Fatou's lemma,
for $y\le -1$,
\[
\E\mu_t(-\iy,y]\le\E\mu_t(-\iy,y+1)\le c_3e^{-c_2(y+1)^2}.
\]
Thus $\E\int_{-\iy}^{-1}\mu_t(-\iy,y]dy<\iy$, which proves the assertion.

Next consider the final assertion.
For $y\in\N$, and $t\in[0,T]$,
if $x^n_i\ge j\ge y\ge X^n_i(t)$ then $\|W_i\|^*_T\ge j-y$,
hence
\[
\mu^n_t(-\iy,y]=\frac{1}{n}\sum_{i,j\ge0}1_{\{X^n_i(t)\le y,x^n_i\in[j,j+1]\}}
\le \mu^n_0[0,y+1]+\frac{1}{n}\sum_{i\ge0,j\ge y}
1_{\{\|W_i\|^*_T\ge j-y\}}1_{\{x^n_i\in[j,j+1]\}}.
\]
Thus, using \eqref{eq12.1}, for some positive constants $c_i$, $i= 1, 2, \ldots, 5$,
\begin{align*}
\E\sup_{t\le T}\mu^n_t(-\iy,y]&\le \E\mu^n_0[0,y+1]
+c_1\sum_{j\ge y}e^{-c_2(j-y)^2}\E \mu^n_0[j,j+1]
\\
&\le c_3(1+y^{m+1})+c_3\sum_{j\ge y}e^{-c_2(j-y)^2}(1+j^m)\\
&\le c_3(1+y^{m+1}) + c_4 \sum_{j\ge 0}e^{-c_2j^2}(1+j^m + y^m) \le c_5(1+ y^{m+1}).
\end{align*}
Next, apply Fatou's lemma to conclude, for some $\tilde c>0$, 
$$
\E\sup_{t\le T}\mu_t(-\iy,y]\le \tilde c(1+y^{\tilde c}),\qquad y\in\N.
$$
By Borel-Cantelli, if $p(y)=\tilde c(1+y^{\tilde c})y^2$,
then $\sup_{t\le T}\mu_t(-\iy,y]>p(y)$ for only finitely many $y\in\N$.
Hence there exists a (random) c such that
$\sup_{t\le T}\mu_t(-\iy,y]\le c(1+y^c)$ for all $y\ge 0$.

\qed

\begin{lemma}\label{lem2n}
Let Assumption \ref{ass1} hold.\\
(a) For every $T\ge 0$ there is an $r>0$ such that
$\lim_{n\to \infty} \sup_{t \in [0,T]}\PP(Y^n_0(t) > r)=0$.
\\
(b) The sequence $\{\beta^n\}_{n\in\N}$ is tight in
$\calM^{(1)}(\R\times\R_+)$.
\end{lemma}

\proof
(a)
By Assumption \ref{ass1}(a) there exists $r_0$ for which
$\PP(\mu^n_0[0,r_0]\le 1)\to0$ as $n\to\iy$.
Note that
$Y^n_0(t)\le n^{-1}\sum_{i<n}X^n_i(t)$. Moreover,
by \eqref{61},
\[
\sum_{i<n}X^n_i(t)\le nt+\sum_{i<n}(x^n_i+W_i(t))
\le nt+nx^n_n+\sum_{i<n}W_i(t).
\]
Hence
	\begin{align*}
		\PP(Y^n_0(t) > 2r_0+t)  &\le 
\PP\Big( nx^n_n+\sum_{i<n} W_i(t) > 2nr_0\Big)\\
		&\le \PP(x^n_n > r_0) + \PP(W(nt)  > nr_0)\\
		&=\PP(\mu^n_0[0,r_0]\le 1)+\PP(W(nt)>nr_0).
	\end{align*}
    Thus we have
    $$ \sup_{t \in [0,T]} \PP(Y^n_0(t) > 2r_0+T) \le \PP(\mu^n_0[0,r_0]\le 1)+\sup_{t \in [0,T]}\PP(W(nt)>nr_0).$$
Sending $n\to \infty$, the result follows with $r=2r_0+T$.

(b)
Fix $T$.
We will have tightness of
$\beta^n$ if we show that for every $\eps$ there are
 $r, \tilde r \in \R$ such that
\[
\limsup_n\PP(\beta^n((-\iy,r]\times[0,T])>\eps)<\eps, \;\; \limsup_n\PP(\beta^n((\tilde r, \infty)\times[0,T])>\eps)<\eps.
\]
For the second statement, note that 
\begin{align*}
\PP(\beta^n((\tilde r, \infty)\times[0,T])>\eps)
&\le \PP \left(\int_0^T 1_{\{Y^n_0(t) >\tilde r\}} dt > \eps\right) \\
&\le \eps^{-1} \int_0^T \PP(Y^n_0(t) > \tilde r) dt
\le \eps^{-1} T \sup_{t\in [0,T]} \PP(Y^n_0(t) >\tilde r).
\end{align*}
The second statement now follows from (a).  Consider now the first statement.
Let, for $j \in \N$,
\[
\beta^{n,j}(dx,dt)=\sum_{i:i\le j}
\del_{X^n_i(t)}(dx)1_{\{X^n_i(t)=Y^n_0(t)\}}.
\]
With $R=(-\iy,r]\times[0,T]$, we have a.s.\ that
$\beta^{n,j}(R)\uparrow\beta^n(R)$ as $j\to\iy$ and
$\PP(\beta^n(R)>\eps)=\lim_j\PP(\beta^{n,j}(R)>\eps)$.
Hence it suffices to show that
\begin{equation}\label{eq:419}
\limsup_n\lim_j\PP(\beta^{n,j}(R)>\eps)<\eps.
\end{equation}
Toward this goal, recall by Lemma \ref{lem1}(c)
that given $\eps \in (0,1)$ there is $r$ such that
\begin{equation}\label{b1}
\sup_n\PP(\sup_{t\le T}\mu^n_t(-\iy,r]>\eps)<\eps.
\end{equation}
Recall that
\[
X_i^n(t)=x^n_i+W_i(t)+n\int_0^t1_{X^n_i(s)=Y^n_0(s)}ds.
\]
Fix $r \in \R$.
Let $\ph = \ph(\cdot;r)\in C^\iy(\R)$, be s.t.\ $\ph=3$ on $(-\iy,r)$,
$\ph=0$ on $(r+3,\iy)$, monotone decreasing
in $[r,r+3]$, with $\ph'=-1$ in $[r+1,r+2]$.
Then for some $c>0$, $|\ph'|<c$, and $|\ph''|<c$. If $r<-3$, $\ph(x^n_i)=0$ for all $i \in \N_0$
and thus
\[
\ph(X^n_i(t))=\int_0^t\ph'(X^n_i(s))(dW_i(s)+n1_{X^n_i(s)=Y^n_0(s)}ds)+\frac12\int_0^t\ph''(X^n_i(s))ds.
\]
Since
\[
-\sum_{i:i\le j}\int_0^t\ph'(X^n_i(s))1_{X^n_i(s)=Y^n_0(s)}ds
=-\int_{\R\times[0,t]}\ph'(x)\beta^{n,j}(dx,ds)
\ge\beta^{n,j}([r+1,r+2]\times[0,t]),
\]
and $|\ph''(x)|\le c1_{[r,r+3]}(x)$, we have
\[
\beta^{n,j}([r+1,r+2]\times[0,t])
\le \frac{1}{n}\sum_{i:i\le j}\int_0^t\ph'(X^n_i(s))dW_i(s)+\frac{c}{2}\int_0^t\mu^n_s[r,r+3]ds.
\]
Summing over $r=r_0, r_0-1,r_0-2,\ldots$, for some $r_0<-3$,
\[
\beta^{n,j}((-\iy,r_0+2]\times[0,t]) \le M^{n,j}(t)+J^n(t),
\]
\[
M^{n,j}(t)= \frac{1}{n}\sum_{i:i\le j}\int_0^t\psi(X^n_i(s))dW_i(s),
\qquad
J^n(t)=\frac{3c}{2}\int_0^t\mu^n_s(-\iy,r_0+3]ds,
\]
\[
\psi=\ph'(\cdot;r_0)+\ph'(\cdot;r_0-1)+\cdots,\qquad
|\psi(x)|\le 3c1_{x\le r_0+3}.
\]
Fix $\eps>0$.  Then, 
by \eqref{b1},  there is $r_0< -3$ s.t. $\PP(J^n(T)>\eps)<\eps$
for all $n$ and $$\sup_n\PP(\sup_{t\le T}\mu^n_t(-\iy,r_0+3]\ge 1)<\eps.$$
Also, $M^{n,j}$ is a martingale on the
filtration generated by $(X^n_i,W_i)_{i\in\N_0}$, and
\[
[M^{n,j}](t)=n^{-2}\sum_{i:i\le j}\int_0^t\psi(X_i(s))^2ds
\le n^{-1}(3c)^2\int_0^t\mu^n_s(-\iy,r_0+3]ds.
\]
Let $\tau=\tau^n=\inf\{t:\mu^n_t(-\iy,r_0+3]\ge1\}$,
and note that it is a stopping time on the aforementioned
filtration. Moreover, by our choice of $r_0$, $\PP(\tau\le T)<\eps$. Hence
\begin{align*}
\PP(M^{n,j}(T)>\eta)&\le\PP(\tau\le T)+\PP(\|M^{n,j}\|^*_{T\w\tau}>\eta)
\\
&
\le\eps+4\eta^{-2}\E\{[M^{n,j}](T\w\tau)\}
\le\eps+4(3c)^2\eta^{-2}n^{-1}T.
\end{align*}
Sending $j\to \infty$ and then $n \to \infty$, and combining this with the above estimate on $J^n(T)$, we have \eqref{eq:419}, completing the proof of the lemma.
\qed

Combining Lemma \ref{lem1} and Lemma \ref{lem2n} we have the following corollary.
\begin{corollary}\label{cor1}
Under Assumption \ref{ass1},
$(\mu^n,\beta^n)$ is a tight sequence in
$C(\R_+,\calM_*)\times\calM^{(1)}(\R\times\R_+)$.
\end{corollary}

\section{Particle system asymptotics}\label{sec3}

The main goal of this section is to prove
that limits are supported on solutions to \eqref{MFBP},
stated as follows.

\begin{lemma}\label{lem2}
Let Assumption \ref{ass1} hold.
Let $(\mu,\beta)$ be a subsequential weak limit of $(\mu^n,\beta^n)$ and
let $v(\cdot,t)=\bF^{\mu_t}$, $t>0$. Then, a.s.,
\\
(a) The function $v$ satisfies, for every $t>0$,
\begin{equation}\label{72}
\lim_{x\to-\iy}v(x,t)=0,
\qquad
\sup_{s\le t}v(x,t)\le c(1+x_+^c),
\end{equation}
for some $c=c(t)$. Moreover, $(v,\beta)$ satisfy \eqref{03},
with data $v_0= \bF^{\mu_0}$,
for all $\ph\in C_c^\iy(\R\times\R_+)$.
\\
(b) $v\in C^\pol_{\#,\loc}(\R\times(0,\iy))$.
\\
(c)
$(v,\beta)$ is a solution to \eqref{MFBP} with initial data
$v_0$.
\end{lemma}

The following lemma is needed in order to prove Lemma \ref{lem2},
specifically, to deduce (b) from (a). It is also crucially used in Section \ref{sec4}. It provides a representation, known as a mild solution, and based on Duhamel's principle, for functions $v$ satisfying \eqref{03} in terms of $v_0$ and $\beta$. Here, $v$ need not a priori be continuous.

For topological spaces $\XX$ and $\YY$, we  denote the space of Borel maps from $\XX$ to $\YY$ by $\calB(\XX, \YY)$. When $\YY = \R$ we simply write $\calB(\XX)$.

Denote by $\frp_t(\cdot-x)$ the Green function associated with $\frac{1}{2}\pl_{xx}$ (namely, the probability density
function of $x+W_t$), and by
$\{\calS_t\}$ the corresponding semigroup:
\begin{align*}
\calS_tf(x) &= \int_{\R} \frp_t(y-x) f(y) dy
\qquad f\in\calB(\R,\R_+)
\\
\calS_t\mu(x) &= \int_{\R} \frp_t(y-x) \mu(dy)
\hspace{2.5em} \mu \in \calM_+(\R),
\end{align*}
where $\calM_+(\R)$ is the space of nonnegative measures on $\R_+$.

\begin{lemma}\label{lem:mild} 
Let $v_0\in\calV_0^\pol(\R)$ and
$\beta\in\calM^{(1)}(\R\times\R_+)$. Let
$v\in\calB(\R\times(0,\iy),\R_+)$ satisfy \eqref{72}
as well as \eqref{03} for all $\ph\in C_c^\iy(\R\times\R_+)$.
Then for every $\tau\in(0,\iy)$ and a.e.\ $(x,t)\in\R\times(0,\iy)$,
\begin{align}
\label{82}
v(x,t) &= \calS_t v_0(x) - \int_0^t \calS_{t-s} \beta_s (x) ds,
\\
\label{83}
v(x,t+\tau) &= \calS_t v_{\tau}(x) - \int_{\tau}^{t+\tau} \calS_{t+\tau-s} \beta_s (x) ds,
\end{align}
where $v_{\tau}(\cdot) = v(\cdot, \tau)$.
Moreover, the r.h.s.\ of \eqref{82} belongs to $C(\R\times(0,\iy))$.
\end{lemma}

\begin{remark}
For $(v,\beta)$ as in Lemma \ref{lem2}, the formulas
\eqref{82}--\eqref{83} are satisfied {\it everywhere}
in $\R\times(0,t)$, as clearly follows
from the a.e.\ statement of Lemma \ref{lem:mild}, the continuity in $(t,x)$ of the right sides of \eqref{82}--\eqref{83},
and the continuity stated in Lemma \ref{lem2}(b).
\end{remark}

\noi{\bf Proof of Lemma \ref{lem:mild}.}
We focus on proving \eqref{82}, the proof of \eqref{83}
being similar.
Given $v_0$ and $\beta$, let the r.h.s.\ of \eqref{82}
be denoted by $\bar v(x,t)$. 
Using dominated convergence, it is easily checked that $\bar v$ satisfies the first statement in \eqref{72} and using the growth condition on $v_0$ it follows that the second statement in \eqref{72} holds for $\bar v$ as well.

An elementary calculation
via integration by parts shows that, for every test function
$\ph$ as in Definition \ref{def0}, \eqref{03} holds for $\bar v$.
Therefore \eqref{82} will follow once uniqueness is
established. The precise uniqueness statement required is:

\noi
{\it Given $v_0$ and $\beta$ as in the statement of the lemma, let $u$ and $v$ be members of
$\calB(\R\times(0,\iy),\R_+)$ satisfying \eqref{72}
as well as \eqref{03}.
Then $u=v$ a.e.\ in $\R\times(0,\iy)$.}

To prove this statement, let
$w=u-v$. Then for $\ph\in C^\iy_c(\R\times\R_+)$,
\[
\int_{\R\times\R_+}(\ph_t+\half\ph_{xx})wdxdt=0.
\]
Extend $w$ to $\R^2$ by setting $w=0$ for $t\le0$.
Then for $\ph\in C^\iy_c(\R^2)$,
\[
\int_{\R^2}(\ph_t+\half\ph_{xx})wdxdt=0.
\]
Let $\eta^{(\eps)}=\eta^{(\eps)}(x,t)$ be the standard (symmetric) mollifier on $\R^2$,
supported on $B_\eps(0) = \{y \in \R^2: \|y\| \le \eps\}$. Let $\psi\in C^\iy_c(\R^2)$.
Let $\ph^{(\eps)}=\eta^{(\eps)}*\psi$, where the convolution is in both
variables. Then
$\ph^{(\eps)}_t=\eta^{(\eps)} *\psi_t$ and $\ph^{(\eps)}_{xx}=\eta^{(\eps)} *\psi_{xx}$,
hence
\[
\int_{\R^2}\eta^{(\eps)}*(\psi_t+\half\psi_{xx})wdxdt=0.
\]
Let $w^{(\eps)}=\eta^{(\eps)}*w$. Then
\[
\int_{\R^2}(\psi_t+\half\psi_{xx})w^{(\eps)}dxdt=0.
\]
Since $w^{(\eps)}\in C^\iy(\R^2)$, integration by parts shows that it satisfies classically
$w_t-\frac{1}{2}w_{xx}=0$ in all of $\R^2$. In particular, it satisfies this equation
in $\R\times[-\eps,\iy)$ with initial condition $w^{(\eps)}(x,-\eps)=0$ for all
$x$.
By uniqueness of classical solutions to the heat equation,
with growth condition in $x$ (for example, \cite[Theorem 7 p.\ 58]{evans97},
we have $w^{(\eps)}=0$ on $\R\times[-\eps,\iy)$.
But $w^{(\eps)}\to w$ a.e. hence $w=0$ a.e., proving the claim.

Next we prove that $\bar v\in C(\R\times(0,\iy))$.
Let $v_1(x,t)=\calS_tv_0(x)$ and $v_2(x,t)=\int_0^t\calS_{t-s}\beta_s(x)ds$. It is well-known that $v_1$ is in $C^\iy$
\cite[Lemma 3.2, p.\ 68]{godlewski91}. However, we include a proof of continuity, as it is simple. Let $\R\times(0,\iy)\ni(x_n,t_n)\to(x,t)\in\R\times(0,\iy)$. Then $\frp_{t_n}(y-x_n)\to\frp_{t}(y-x)$, and since
$(x_n,t_n)$ are contained in a compact subset of $\R\times(0,\iy)$,
\[
\frp_{t_n}(y-x_n)v_0(y)= ct_n^{-1/2}e^{-(y-x_n)^2/t_n}v_0(y)\le c_1e^{-c_2y^2}(1+y_+^c).
\]
Hence $\calS_{t_n}v_0(x_n)\to\calS_tv_0(x)$ by dominated convergence.

As for $v_2$, for $0<\eps<\liminf t_n/2$, and $n$ large,
\[
\int_{t_n-\eps}^{t_n}\int_\R\frp_{t_n-s}(y-x_n)\beta_s(dy)ds\le\int_{t_n-\eps}^{t_n}\int_\R
(t_n-s)^{-1/2}\beta_s(dy)ds\le\eps^{1/2}.
\]
A similar estimate holds when $(x_n,t_n)$ is replaced by $(x,t)$. Hence it suffices to show that      
for each $\eps$ as above,
\[
\int_0^{t_n-\eps}\int_\R\frp_{t_n-s}(y-x_n)\beta_s(dy)ds \to
\int_0^{t-\eps}\int_\R\frp_{t-s}(y-x)\beta_s(dy)ds.
\]
But $\frp_{t_n-s}(y-x_n)\le\eps^{-1/2}$,
hence the above convergence holds by dominated
convergence, recalling that $\beta_s$ are probability measures. This completes the proof that $\bar v$ is continuous.
\qed

\noi{\bf Proof of Lemma \ref{lem2}.}
(a) Properties \eqref{72} are immediate from Lemma \ref{lem1}. Next, to verify \eqref{03},
let $\ph\in C^\iy_c(\R\times\R_+)$. Then for some $r$ and $T$,
$\ph$ vanishes on $\{t> T\}$ and on $\{x> r\}$.
Let $\psi(x,t)=-\int_x^\iy\ph(y,t)dy$. Then $\psi$ also vanishes
on the above domains, and is bounded on $\R\times\R_+$.
For $i\in\N_0$, by Ito's lemma,
\begin{align*}
0=\psi(X^n_i(T),T)&=\psi(x^n_i,0)+\int_0^T(\pl_t\psi+\half\pl_{xx}\psi)(X^n_i(t),t)dt
\\
&\quad+n\int_0^T1_{\{X^n_i(t)=Y^n_0(t)\}}\pl_x\psi(X^n_i(t),t)dt+M^n_i(T),
\end{align*}
where $M^n_i(t)=\int_0^t\pl_x\psi(X^n_i(s),s)dW_i(s)$.
Summing over $i$ and dividing by $n$,
\[
0=\lan\psi(\cdot,0),\mu^n_0\ran
+\int_0^T\lan \pl_t\psi+\half\pl_{xx}\psi,\mu^n_t\ran dt
+\int_{\R\times[0,T]}\pl_x\psi(x,t)\beta^n(dx,dt)
+\bar M^n(T)
\]
where $\bar M^n(t)=n^{-1}\sum_{i\ge0}M^n_i(t)$. Observe that
\[
[\bar M^n]_T\le cT n^{-2}
\sum_{i\ge0}1_{\{\min_{[0,T]}X^n_i<r \}}.
\]
Taking expectation and using Lemma \ref{lem1}(b) shows that
$\|\bar M^n\|^*_T\to0$ in probability.
Consider a subsequential limit $(\mu,\beta)$. Then it satisfies
\[
0=\lan\psi(\cdot,0),\mu_0\ran
+\int_0^T\lan \pl_t\psi+\half\pl_{xx}\psi,\mu_t\ran dt
+\int_{\R\times[0,T]}\pl_x\psi(x,t)\beta(dx,dt).
\]
Here the convergence of 
$\int_0^T\lan \pl_t\psi+\half\pl_{xx}\psi,\mu^n_t\ran dt$
to 
$\int_0^T\lan \pl_t\psi+\half\pl_{xx}\psi,\mu_t\ran dt$  follows on recalling that $\mu^n \to \mu$ in $C(\R_+, \mathcal{M}_*)$ (along the subsequence, in probability), and  $f = \pl_t\psi+\half\pl_{xx}\psi$
satisfies $\sup_{t\in [0,T]} \|f(t, \cdot)\|_{\BL} < \infty$ and, for some $r<\infty$, $f(t, \cdot)$ has support contained in $(-\infty, r]$, for every $t\in [0,T]$.
Recall $v_0=\bF^{\mu_0}$ and $v(\cdot,t)=\bF^{\mu_t}$.
Then integration by parts, using that
$\psi(x,t)$ vanishes for $x>r$ and $v(x,t)$ vanishes in the limit $x\to-\iy$,
gives
\[
0=-(\ph(\cdot,0),v_0)-\int_0^T(\pl_t\ph+\half\pl_{xx}\ph,v(\cdot,t))dt
+\int_{\R\times[0,T]}\ph(x,t)\beta(dx,dt).
\]
$T$ can now be replaced by $\iy$, and \eqref{03} is established.

(b) First we show that $v\in C(\R\times(0,\iy))$.
From the a.e.\ equality in \eqref{82}, we can find $\Lambda \in \calB(0, \infty)$ such that the Lebesgue measure of $\Lambda^c$ is $0$ and for all $t \in \Lambda$, \eqref{82} holds for a.e.\ $x \in \R$.
Denote the r.h.s.\ of \eqref{82} by $\bar v(x,t)$ and by $\bar\mu_t$, $t\ge0$, the measure $\bar\mu_t(-\iy,x]=\bar v(x,t)$. Then
given a test function $\ph\in C_c(\R,\R)$, $\lan\ph,\mu_t\ran=\lan\ph,\bar\mu_t\ran$ for a.e.\ $t$. But $t\mapsto\mu_t$ is continuous in $\calM_*$ (cf.\ Corollary \ref{cor1}), and clearly so is $t\mapsto\bar\mu_t$. This gives $\mu_t=\bar\mu_t$ for all $t$. Hence $v=\bar v$ holds in $\R\times(0,\iy)$, and therefore by Lemma \ref{lem:mild}, $v\in C(\R\times(0,\iy))$.

In view of Lemma \ref{lem1}(e), in order to prove
that $v\in C^\pol_{\#,\loc}(\R\times(0,\iy))$,
it only remains to show that $\lim_{x\to\iy}v(x,t)=\iy$, $t>0$.
To this end, use \eqref{82} and the estimates
\[
\int_0^t\calS_{t-s}\beta_s(x)ds=\int_0^t\int_\R\frp_{t-s}(y-x)\beta_s(dy)ds
\le \int_0^tc(t-s)^{-1/2}ds=ct^{1/2},
\]
\[
\int_\R\frp_t(y-x)v_0(y)dy\ge v_0(x)\int_x^\iy\frp_t(y-x)dy
=\frac{v_0(x)}{2},
\]
where the fact that $v_0$ is nondecreasing is used.
The result follows from $\lim_{x\to\iy}v_0(x)=\iy$.

(c)
It remains to show that $(v,\beta)$ satisfy (\ref{MFBP}.iii).
Invoking Skorohod's representation, we may assume without loss of generality that the subsequential convergence $(\mu^n,\beta^n)\to(\mu,\beta)$ holds a.s.
Denote $v^n(\cdot,t)=\bF^{\mu^n_t}$ and
$u^n(x,t)=v^n(x-,t)$. By the definition of $\beta^n$,
it is supported on the set $\{(x,t):x=Y^n_0(t)\}$, whereas
for every $t$,
\[
u^n(Y^n_0(t),t)=\mu^n_t(-\iy,Y^n_0(t))=0.
\]
This gives
\[
\int_{\R\times\R_+}u^n(x,t)\beta^n(dx,dt)=0.
\]
Thus the claim will follow once it is established that
\begin{equation}\label{81}
\int_{\R\times\R_+}vd\beta
\le\liminf_n\int_{\R\times\R_+}u^nd\beta^n.
\end{equation}

Fix $r,T\in(0,\iy)$ and let $G=(-r,r)\times(0,T)$.
Given $\eps \in (0,1)$ and $x \in (-r,r)$, consider a function $f^{\eps,x}: \R \to [0,1]$, with Lipschitz constant bounded by $2\eps^{-1}$, such that 
$f^{\eps,x}(y) = 1$ for $y \in (-\iy, x-\eps)$ and $0$ for $y \ge x-\eps/2$.
Then, by uniform convergence on compacts of $\mu^n$ to $\mu$, we can find an $n_0 \in \N$ such that for all $n\ge n_0$
$$\sup_{t\in [0,T]} \sup_{x \in(-r,r)} | \langle f^{\eps, x}, \mu^n_t\rangle - \langle f^{\eps, x}, \mu_t\rangle|< \eps.$$
Here we are using the fact that the collection $\{ \frac{\eps}{2}f^{\eps, x}, x \in (-r,r)\}$ has $\text{BL}$ norm bounded by $1$ and support contained in $(-\infty, r+1]$.
Thus, for all $(x,t) \in G$ and $n \ge n_0$,
$$
v(x-\eps, t) \le \langle f^{\eps,x}, \mu_t\rangle \le \langle f^{\eps,x}, \mu^n_t\rangle +\eps \le u^{n}(x,t) + \eps.$$
Since $v$ is continuous and $\beta^n\to\beta$ weakly, we have, with $v_{\eps}(\cdot) = v(\cdot-\eps)$
\[
\int_G v_{\eps} d\beta \le\liminf_n\int_Gv_{\eps} d\beta^n
\le\liminf_n \int_G u^nd\beta^n + \eps T.
\]
Now, as $\eps\to0+$, $v_{\eps}\to v$ pointwise. Hence by bounded
convergence, $\int_G v_{\eps}d\beta\to\int_Gvd\beta$.
Since $r$ and $T$ are arbitrary, the inequality in \eqref{81} follows.
This proves (\ref{MFBP}.iii) and completes the proof
of the result.
\qed

\section{PDE uniqueness}\label{sec4}

In this section we prove the uniqueness statement in Theorem \ref{th1}(a).

\begin{theorem}
\label{thm:uniq}
Let $v_0 \in \calV^\pol_0(\R)$ and let $(v^i, \beta^i) \in C^\pol_{\#,\loc}(\R\times(0,\iy))\times\calM^{(1)}(\R\times\R_+)$, $i=1,2$, be two solutions of \eqref{MFBP}  with initial condition $v_0$. Then
$(v^1, \beta^1) = (v^2, \beta^2)$.
\end{theorem}

The technique uses ideas from \cite{ata25},
where uniqueness was proved for a related equation
for solutions in the class $L_{1,\loc}(\R_+,L_q(\R))$,
$q\in(1,\iy)$. The function $v$ of interest in this paper
is clearly not in this class. This necessitates various
significant adaptations of the argument.

In what follows we use some notation special to this 
section. For a function $f:\R\to\R$, $\int_a^bf(x)dx$
is abbreviated as $f\cc{a,b}$.
Let $\calB^\pol_\#(\R)$
denote the collection of $v\in\calB(\R,\R_+)$ satisfying
\[
v(x)\le c(1+x_+^c),\qquad x\in\R, \qquad
\lim_{x\to\iy}v(x)=\iy, \qquad
v|_{(-\iy,0)}\in L^1(-\iy,0).
\]
Clearly, intersecting
$\calB^\pol_\#(\R)$ with the class of continuous
nondecreasing functions gives $C^\pol_\#(\R)$ and intersecting with the class of right continuous non decreasing functions gives $\calV^\pol_\#(\R)$.
The reason for introducing $\calB^\pol_\#(\R)$ is that the tools to be used here require us to work in a class of functions that need not be monotone.
Throughout what follows, $\calB^\pol_\#(\R)$ is abbreviated to $\calB^0$.

For $\delta>0$, define the maps 
$\calC_{\delta}: \calB^0 \to \calB^0$  as
\begin{equation*}
\calC_{\delta} v\doteq 
  v(x) \one_{[\gamma^{\delta}_v, \infty)}(x),\qquad
 \text{ where } \gamma^{\delta}_v = \inf\{a \in \R: v[-\infty, a] = \delta\}.
\end{equation*}
Note that these are well-defined owing to the condition $v|_{(-\iy,0)}\in L^1$.
For $v_1, v_2 \in \calB^0$ and $\ell \ge 0$ we denote $v_1\LE v_2 \mod \ell$ if
$$v_2\cc{-\infty,r} \le v_1\cc{-\infty,r}  + \ell \mbox{ for all } r \in \R.$$
We abbreviate $v_1\LE v_2 \mod 0$ to $v_1\LE v_2$. We note that, for any $v \in \calB^0$ and $r \in \R$, $v[r, \infty] = \infty$.

\begin{lemma} 
\label{lem:simpest}
Let $\delta>0$ and $\ell\ge0$. Then we have the following.
\begin{enumerate}[(a)]
\item For $v \in \calV^\pol_\#(\R)$ (resp. $\calB^0$), $\calC_{\delta} v \in \calV^\pol_\#(\R)$ (resp. $\calB^0$) and $\calS_{\delta} v \in \calV^\pol_\#(\R)$ (resp. $\calB^0$).
\item For $u, v \in \calB^0$ with $u\LE v \mod \ell$,
$\calC_{\delta} u \LE \calC_{\del} v \mod \ell$.
\item  For $u, v \in \calB^0$ with $u\LE v \mod \ell$,
$\calS_{\delta} u \LE \calS_{\delta} v \mod \ell$.
\end{enumerate}

\end{lemma}

\proof
(a) This statement regarding $\calC_{\delta}$ is immediate from the definition.
Also,
$$\calS_{\delta}v(x) = \int_{\R} v(y) \frp_{\del}(x-y) dy = \int_{\R} v(x+z) \frp_{\del}(-z) dz.$$
Therefore the polynomial growth (resp. monotonicity) of $v$ implies that of
$\calS_\del v$. Using dominated convergence it follows that $\calS_{\delta}v(x) \to \infty$ as $x \to \infty$. Moreover,
\begin{align*}
\calS_{\delta}v\cc{-\infty,0} &= \int_{\R} \frp_{\del}(-z) \int_{-\infty}^{z} v(x) dx\, dz\\ 
&\le \int_{\R}\frp_{\del}(-z) \int_{-\infty}^{0} v(x) dx\, dz + \int_{\R} \frp_{\del}(-z) \int_{0}^{z^+} v(x) dx\, dz<\iy,
\end{align*}
owing to $v\cc{-\iy,0}<\iy$ and the polynomial growth
of $v$.

(b) 
We have, for every $r \in \R$
$v\cc{-\infty,r} \le u\cc{-\infty,r} + \ell.$
Note that $\calC_{\del} u = u \one_{[\gamma_u^{\del}, \infty)}$ and
$\calC_{\del} v = v \one_{[\gamma_v^{\del}, \infty)}$. 
Then the inequality
$\calC_{\del}v\cc{-\infty,r} \le  \calC_{\del}u\cc{-\infty,r} + \ell$
is immediate for $r< \gamma_v^{\del}$.
Also, for $r\ge \gamma_v^{\del}$,
\begin{align*}
 \calC_{\del} v\cc{-\infty,r} &= v\cc{-\infty,r} - v\cc{-\infty,\gamma_v^{\del}}
= v\cc{-\infty,r} - u\cc{-\infty,\gamma_u^{\del}} \\ 
&\le
 u\cc{-\infty,r} +\ell - u\cc{-\infty,\gamma_u^{\del}} 
 \le   \calC_{\del} u\cc{-\infty,r} +\ell,
\end{align*}
where the last inequality uses the observation that if $r < \gamma_u^{\del}$ then $\calC_{\del} u\cc{-\infty,r} =0$ and $u\cc{-\infty,r}  - u\cc{-\infty,\gamma_u^{\del}} \le 0$.
This proves  (b).

 (c) Note that, for $u,v$ as in the statement,
\begin{align*}
\calS_{\del} v\cc{-\iy,r} &= \int_{-\infty}^r \int_{\R} v(x+z) \frp_{\del}(-z) dz \, dx\\ 
&= \int_{\R} \frp_{\del}(-z) \int_{-\infty}^{r+z} v(x) dx\, dz
\le \int_{\R} \frp_{\del}(-z) [\int_{-\infty}^{r+z} u(x) dx + \ell] dz\\ 
&= \int_{\R} \frp_{\del}(-z) \int_{-\infty}^{r+z} u(x) dx \, dz + \ell =
\calS_{\del} u\cc{-\iy,r} + \ell.
\end{align*}
\qed

Given $\Del>0$ and $\delta\ge 0$, denote $\hat\delta = (\Del, \delta)$ and define
the operator $\calC_{\hat\delta}:\calB^0 \to \calB^0$ as follows.
$$\calC_{\hat\delta}v(x) = v(x)\one_{(-\infty, \gamma^{\Del}_v)}(x) + v(x) \one_{[\gamma^{\Del+\delta}_v, \infty)}(x), \; x \in \R, \; v \in \calB^0.$$
Occasionally we write $\calC_{\Delta,\delta}$ instead of $\calC_{\hat\delta}$.

\begin{lemma}
\label{lem:q1}
Let $\ell\ge0$, $\delta, \Del>0$ and $u, v \in \calB^0$ such that $u\LE v \mod \ell$. Then
$\calC_{\hat\delta}u \LE \calC_{\hat\delta}v \mod \ell$.
\end{lemma}

\proof
Fix $r \in \R$. It is required to show that
\begin{equation}
 \calC_{\hat\delta}v\cc{-\infty,r}\le \calC_{\hat\delta}u\cc{-\infty,r} + \ell.
\end{equation}
We consider various cases.
\\
Case 1: $r \in [\gamma_{v}^{\Del+\delta}, \infty)$.
Note that, with $\calD_{\hat\delta}u = u- \calC_{\hat\delta}u$
\begin{align*}
\calC_{\hat\delta}v\cc{-\iy,r} &= v\cc{-\iy,r} - 
\calD_{\hat\delta}v\cc{-\iy,r} = v\cc{-\iy,r}
-
\delta \\ 
&\le u\cc{-\iy,r} - \delta +\ell = \calC_{\hat\delta}u\cc{-\iy,r}
+ \calD_{\hat\delta}u\cc{-\iy,r} - \delta +\ell\\ 
&\le \calC_{\hat\delta}u\cc{-\iy,r} + \delta - \delta +\ell =
\calC_{\hat\delta}u\cc{-\iy,r} +\ell.
\end{align*}
Case 2: $r \in (-\infty, \gamma_{v}^{\Del}]$. 
Consider the subcase when one also has $r \in (-\infty, \gamma_{u}^{\Del}]$. In this case
\begin{align*}
\calC_{\hat\delta}v\cc{-\iy,r} &=  v\cc{-\iy,r} 
\le u\cc{-\iy,r} +\ell = \calC_{\hat\delta}u\cc{-\iy,r} +\ell.
\end{align*}
Now consider the subcase $r \in (\gamma_{u}^{\Del}, \infty)$. Then 
\begin{align*}
\calC_{\hat\delta}v\cc{-\iy,r} &\le v\cc{-\infty,\gamma_{v}^{\Del}}
= u\cc{-\infty,\gamma_{u}^{\Del}}
= \calC_{\hat\delta}u\cc{-\infty,\gamma_{u}^{\Del}}
\le \calC_{\hat\delta}u\cc{-\iy,r}
\le \calC_{\hat\delta}u\cc{-\iy,r}+\ell.
\end{align*}
Case 3: $r \in (\gamma_{v}^{\Del}, \gamma_{v}^{\Del+\delta})$. Consider first the subcase where $r \in (-\infty, \gamma_{u}^{\Del}]$. Then
\begin{align*}
\calC_{\hat\delta}v\cc{-\iy,r} &=
v\cc{-\infty,\gamma_{v}^{\Del}} 
\le u \cc{-\infty,\gamma_{v}^{\Del}} +\ell
\le u\cc{-\infty,r} + \ell
= \calC_{\hat\delta}u\cc{-\iy,r} +\ell.
\end{align*}
Now consider the subcase $r \in [\gamma_{u}^{\Del}, \infty)$. Then
\begin{align*}
\calC_{\hat\delta}v\cc{-\iy,r} &= v\cc{-\infty,\gamma_{v}^{\Del}}
= u \cc{-\infty,\gamma_{u}^{\Del}}
= \calC_{\hat\delta}u \cc{-\infty,\gamma_{u}^{\Del}}
\le \calC_{\hat\delta}u\cc{-\iy,r}
\le \calC_{\hat\delta}u\cc{-\iy,r} + \ell.
\end{align*}
This completes the proof of the lemma. \qed

\begin{lemma}
If $0<\hat \Del \le \Del$ and $v \in \calB^0$, then
$\calC_{\Del, \delta}v \LE \calC_{\hat\Del, \delta}v$.
\end{lemma}

\proof
Fix $r \in \R$. We need to show that
\begin{equation}\label{84}
\calD_{\Del, \delta}v\cc{-\iy,r} \le \calD_{\hat \Del, \delta}v\cc{-\iy,r}.
\end{equation}
Note that $\gamma_{v}^{\hat \Del} \le \gamma_{v}^{\Del} \le \gamma_{v}^{\Del+\delta}$.
We consider the following cases.
\begin{itemize}
\item $r \le \gamma_{v}^{\Del}$. In this case the left-hand side of \eqref{84} is $0$ and so the inequality holds.
\item $r \ge \gamma_{v}^{\Del+\delta}$. Then
both sides of \eqref{84} are equal to $\del$ and once more \eqref{84} holds.
\item $r \in (\gamma_{v}^{\Del}, \gamma_{v}^{\Del+\delta}]$. Consider first the subcase that
$r \le \gamma_{v}^{\hat \Del+\delta}$. Then
\begin{align*}
\calD_{\Del, \delta}v\cc{-\iy,r} &= v\cc{\gamma_{v}^{\Del},r}
\le v\cc{\gamma_{v}^{\hat \Del},r} 
= v\cc{\gamma_{v}^{\hat \Del},r\wedge \gamma_{v}^{\hat \Del+\delta}} 
= \calD_{\hat\Del, \delta}v\cc{-\iy,r}.
\end{align*}
Now consider the subcase $r > \gamma_{v}^{\hat \Del+\delta}$.
Then
$\calD_{\Del, \delta}v\cc{-\iy,r}\le \delta
= \calD_{\hat \Del, \delta}v\cc{-\iy,r}$.
\end{itemize}
This completes the proof of the lemma.
\qed

\begin{lemma}\label{lem:q7}
Suppose $u, v \in \calB^0$ and $\Del, \Del' >0$. If $\Del' \ge \Del$ and $u\LE v \mod \Del'$, then $\calC_{\delta}u \LE \calC_{\Del, \delta}v \mod \Del'$.
\end{lemma}

\proof
We need to show that for every $r\in \R$
\begin{equation}
\calC_{\Del, \delta}v\cc{-\iy,r} \le \calC_{\delta} u\cc{-\iy,r} + \Del'.
\end{equation}
We consider two cases.
\\
Case 1: $r> \gamma_{v}^{\Del +\delta}$. Then
\begin{align*}
\calC_{\Del, \delta}v\cc{-\iy,r} &= v\cc{-\iy,r} - \delta
\le u\cc{-\iy,r} + \Del'- \delta\\
& \le \calC_{\delta}u\cc{-\iy,r} +\delta +\Del'-\delta
= \calC_{\delta}u\cc{-\iy,r} + \Del'.
\end{align*}
Case 2: $r\le  \gamma_{v}^{\Del +\delta}$. In this case
\begin{align*}
\calC_{\Del, \delta}v\cc{-\iy,r} &\le
v\cc{-\infty,\gamma_v^{\Del}}
= \Del \le \Del' \le \calC_{\delta}u\cc{-\iy,r} +\Del'.
\end{align*}
\qed 

Denote $v^{(\del,+)}_0=v_0$ and recursively,
$$
v^{(\del,+)}_{n\del}= \calC_{\delta}\calS_{\delta}v^{(\del,+)}_{(n-1)\del}, \qquad n\in\N.
$$
In what follows, for $n \in \N$, we write $v(\cdot, n\delta)$ as $v_{n\delta}(\cdot)$.
\begin{lemma}\label{lem:q9}
Fix $\delta>0$. Let $(v,\beta)\in C^\pol_{\#,\loc}(\R\times(0,\iy))\times\calM^{(1)}(\R\times\R_+)$
be a solution to \eqref{MFBP} for some $v_0 \in \calV_0^\pol(\R)$.
Then
$$v_{n\delta} \LE v_{n\delta}^{(\delta, +)},\qquad
n\in\N_0.$$
\end{lemma}

\proof
The proof proceeds by induction. Note that $v(\cdot, 0)\doteq v(\cdot, 0+) =v_0(\cdot)$ a.e.\ (this follows, e.g., by taking $t\to0$ in \eqref{82}).
For $n=0$ the claim holds by definition. Suppose the result is true for $(n-1)$ with $n\in \N$.
Let
$$f= v_{(n-1)\delta}, \; g= v_{(n-1)\delta}^{(\delta, +)}, \; \calC = \calC_{\delta}, \; \calS = \calS_{\delta}$$
and assume $f\LE g$.
Note that $v_{n\delta}^{(\delta, +)} = \calC\calS g$ and, from Lemma \ref{lem:mild}, $v_{n\delta} = \calS f -h$, where
$$h(y) \doteq \int_{\R \times [(n-1)\delta, n\delta]} \frp_{n\delta-s}(y-x) \beta(dx\, ds).$$
Thus it suffices to show $\calS f -h \LE \calC\calS g$. Let $w= \calS f$. We claim that $w-h \LE \calC w$. Assuming the claim, we have
$$
\calS f -h = w-h \LE \calC w = \calC \calS f \LE \calC \calS g,$$
where the last inequality uses Lemma \ref{lem:simpest}(b) and (c). Thus the result follows by induction.

To prove the claim, it suffices to show that for all $r\in \R$,
$\calC w\cc{-\iy,r} \le (w-h)\cc{-\iy,r}$.
Note the inequality holds trivially if $r < \gam_w^{\delta}$, since $w-h = v_{n\delta}\ge 0$ by definition of a solution.
Also, if $r \ge \gam_v^{\delta}$,
$$\calC w\cc{-\iy,r}  = w\cc{-\iy,r}
- w\cc{-\infty,\gam_w^{\delta}}
= w\cc{-\iy,r} - \delta \le (w-h)\cc{-\iy,r},$$
completeing the proof, where the last inequality uses the fact that
$h\cc{-\iy,r} \le h\cc{-\iy,\iy} = \delta$.
This completes the proof of the lemma.
\qed

Set $v_0^{(\hat \delta, -)} = v_0$, $\ell_{0, \hat \delta} =0$, and for $n\in\N$,
\begin{equation}\label{eq:aa4}
v_{n\delta}^{(\hat \delta,-)} = \calC_{\hat\delta}\calS_{\delta} v_{(n-1)\delta}^{(\hat \delta,-)},\;\;
\ell_{n, \delh} = \ell_{n-1,\delh} + \delta(\one_{\{(n-1)\del < t_0\}} + e^{-\Del^5/\delta}
\one_{\{(n-1)\del \ge t_0\}}).
\end{equation}
Given $\mu \in \calM^{(1)}(\R\times \R_+)$ and $[t_1, t_2] \subset \R_+$, let
$$\rho^*(\mu; [t_1, t_2]) = \sup \mbox{supp}(\mu(\cdot \times [t_1, t_2])).$$

\begin{lemma}\label{lem:q2}
 Let $(v,\beta)$ be as in Lemma \ref{lem:q9}. Given $\delta>0$ and $n\in \N$, let $\rho_{n,\del} = \rho^*(\beta;[(n-1)\del, n\del])$.
Then, for $n\ge 2$,
$$\rho_{n,\del} \le \bb(n,\del, v_{(n-1)\del}) \doteq \gamma_{v_{(n-1)\del}}^{4\del}.$$
\end{lemma}

\proof
Arguing by contradiction, assume that, for some $n \in \N$,
$$\rho_{n,\del} >\bb = \bb(n,\del, v_{(n-1)\del})= \gamma_{v_{(n-1)\del}}^{4\del}.$$
Then we must have $\theta \doteq \beta((\bb,\infty) \times ((n-1)\del, n\del]) >0$.
Let $V(x,t) \doteq \int_{-\infty}^x v(y, t) dy$ for $(x,t) \in \R\times \R_+$.
We claim that $V(\bb,t)=0$ for some $t \in ((n-1)\del, n\del]$. Suppose this is false. Then for all $t \in ((n-1)\del, n\del]$,
$V(\bb,t)>0$ and so
\begin{align*}\beta(V>0) &\ge \int_{(\bb,\infty)\times ((n-1)\del, n\del]} \one_{V(x,t)>0} \beta(dx\, dt)\\
&\ge \int_{(\bb,\infty)\times ((n-1)\del, n\del]} \one_{V(\bb,t)>0} \beta(dx\, dt) = \theta>0
\end{align*}
which contradicts the definition of $\beta$ (note $v=0$ a.e. $\beta$ says that $V=0$ a.e. $\beta$).
This proves the claim.

Now fix $t \in ((n-1)\del, n\del]$ such that $V(\bb,t)=0$.
We have from Lemma \ref{lem:mild}
$$v(y,t) = \calS_{t-(n-1)\delta} v(\cdot,(n-1)\delta)(y) - S\ast \beta(y,t,(n-1)\delta),
$$
where we denote $S\ast \beta(y,t,(n-1)\delta) = \int_{\R \times [(n-1)\delta, t]} \frp_{t-s}(y-x) \beta(dx\, ds)$.
Letting, for given $t$,
$h(y) = S\ast \beta(y,t,(n-1)\delta)$, we have
$ h\cc{-\infty,\bb} \le h\cc{-\iy,\iy} \le \delta$.
Also, since $V(\bb,t)=0$,
$$ \calS_{t-(n-1)\delta} v(\cdot,\delta)\cc{-\infty,\bb} = \int_{-\infty}^\bb\int_{\R} \frp_{t-(n-1)\del}(y-x) v(x,\del) dx\, dy = h\cc{-\infty,\bb} \le \del.$$
Finally, because for $x\le\bb$, $\int_{-\iy}^\bb\frp_s(y-x)dy\ge\int_{-\iy}^0\frp_s(y)dy=1/2$, we have
\begin{align*}
\del \ge \int_{-\infty}^\bb\int_{-\infty}^\bb \frp_{t-(n-1)\del}(y-x) v(x,\del) dx\, dy
\ge \frac{1}{2}\int_{-\infty}^\bb v(x,\del) dx = \frac{4\del}{2} = 2\del.
\end{align*}
This is a contradiction and the result follows.
\qed

Recall for $\del, \Del>0$, $\hat \del = (\Del, \del)$, and $t_0>0$, $\ell_{n, \hat \del}$ from
\eqref{eq:aa4}.

\begin{lemma}\label{lem:q5}
Let $(v,\beta)$ be as in Lemma \ref{lem:q9}. Let $0< t_0<T$ be given. Then there is a $\Del_0 >0$ such that for every $\Del \in (0, \Del_0)$
there exists $\delta_0 = \del_0(\Del)>0$ such that for $\del \in (0, \del_0)$ and $n \in \N$ s.t. $n\del \le T$,
$$v_{n\delta}^{(\hat \delta,-)} \LE v_{n\del} \mod \ell_{n, \delh}.$$
Also, for $n\in \N$, $n\del \le T$, 
$$\ell_{n, \delh} \le t_0+\del+ Te^{-\Del^5/\del}.$$
\end{lemma}

\proof
Let $\chi_n = \one_{\{(n-1)\del < t_0\}} + e^{-\Del^5/\delta}$. Note that, for $n \in \N$ such that $n\del \le T$,
$$\ell_{n, \delh} \le \ell_{(n-1), \delh} + \del \chi_n.$$
Thus we have
$$\ell_{n, \delh} \le \del \sum_{i=1}^n \chi_i \le t_0+\del + e^{-\Del^5/\delta} T$$
which proves the second statement in the lemma.

Now consider the first part. 
We begin by making a suitable selection of $\Del_0$ and $\del_0(\Del)$. Because $v$ satisfies a polynomial growth condition locally uniformly, we have, for some $c_{v,T}>0$ and $m \in \N$,
\begin{equation}\label{eq:cvt2}
\sup_{t\in [0,T]} v(x,t) \le c_{v,T}(1+ x_+^m), \qquad x \in \R.
\end{equation}
Fix $\Del_*>0$. Then, there exists $\gamma_* \in (0,\infty)$ such that
\begin{equation}
v(\cdot, s)\cc{-\infty,\gamma_*}  \ge \Del_* \mbox{ for all } s \in [0,T].
\end{equation}
This follows on noting that, for any $a \in \R$, and $s \in [0,T]$,
$$\int_{-\infty}^a \int_0^s \calS_{s-u} \beta_u(x) ds \le s$$
and so using \eqref{82}
$$
v(\cdot,s)\cc{-\infty,a} \ge \frac{1}{2}v_0[-\infty, a] - T.$$

Choose 
$$\Del_0 < \min\left\{ \frac{1}{8}, \frac{1}{12} \frac{1}{c_{v,T} (1+ (\gamma_*+1)^m)}, 3\Del_*\right\}.$$
For each fixed $\Del \in (0, \Del_0)$, we choose $\del_0= \del_0(\Del) \in (0,1)$ that satisfies the following conditions:
\begin{itemize}
\item $\sqrt{2} e^{-\Del^5/\del_0} <1$.
\item $\del_0 < \Del/24$.
\item 
$$
\sqrt{2}c_{v,T}\exp\{-\Del^4/(8\del_0)\}\int_{0}^{\infty} 
 (1+ |x+\gamma_*+1|^m) \exp\{-x^2/8\del_0\} dx \le \Del/3.$$
\end{itemize}
Henceforth $\Del \in (0, \Del_0)$ and for such a $\Del$, $\del \in (0, \del_0(\Del))$, where $\Del_0$ and $\del_0(\Del)$ satisfy the above conditions.

We proceed by induction. Suppose $f \LE g \mod \ell_{n-1, \delh}$, where
$$f = v_{(n-1)\delta}^{(\hat \delta,-)}, \; g = v_{(n-1)\del}.$$
Then, $v_{n\delta}^{(\hat \delta,-)} = \calC\calS f$, and from Lemma \ref{lem:mild},  $v_{n\del} = \calS g-h$, where
$$h(y) = \int_{\R \times [(n-1)\del, n\del]} \frp_{n\del-s}(y-x) \beta(dx\, ds).$$
We need to show that
\begin{equation}\label{eq:aa7}
\calC\calS f \LE \calS g -h \mod \ell_{n, \delh}.
\end{equation}
Note from Lemma \ref{lem:simpest}(c) that
$\calS f \LE \calS g \mod \ell_{n-1, \delh}$.
Suppose now that $(n-1)\del < t_0$. Then, for any $r \in \R$,
\begin{align*}
\calC\calS f\cc{-\iy,r} &\ge \calS f\cc{-\iy,r} - \del
\ge \calS g\cc{-\iy,r} - \ell_{n-1, \delh} - \del \\ 
& = (\calS g-h)\cc{-\iy,r} + h\cc{-\iy,r} - \ell_{n,\delh}
\ge (\calS g-h)\cc{-\iy,r} - \ell_{n,\delh}.
\end{align*}
This gives \eqref{eq:aa7} when $(n-1)\del < t_0$.

Now consider $(n-1)\del \ge t_0$.
From Lemma \ref{lem:simpest} and \ref{lem:q1} and induction hypothesis
$$v_{n\delta}^{(\hat \delta,-)} = \calC\calS f \LE \calC\calS g \mod \ell_{n-1, \delh}.$$
Let $w=\calS g$ and $\eps = e^{-\Del^5/\del} \del$. We claim that
\begin{equation}
\calC w \LE w-h \mod \eps . \label{eq:q1}
\end{equation}
Note that once we have the claim, it follows that
$$v_{n\delta}^{(\hat \delta,-)} = \calC\calS f \LE w-h \mod (\ell_{n-1, \delh} + e^{-\Del^5/\del} \del),$$
since $(\ell_{n-1, \delh} + e^{-\Del^5/\del} \del) = \ell_{n, \delh}$, this proves \eqref{eq:aa7}, and the proof is complete by induction.
Thus it remains to show the claim in \eqref{eq:q1}.

Let $\bb = \gamma_g^{4\del}$. Then, from Lemma \ref{lem:q2}, $\rho_{n,\del} \le \bb$.
Write $h = h_1+h_2$, where
$$h_1(y) = h(y) \one_{y > \bb+\Del^2}, \;  h_2(y) = h(y) \one_{y \le \bb+\Del^2}.$$
Since $\rho_{n,\del} \le \bb$
$$h(y) = \int_{(-\infty, \bb]\times [(n-1)\del, n\del]} \frp_{n\del-s}(y-x) \beta(dx\, ds).$$
Note that, with our choice of $\del$ and $\Del$, 
\begin{align*}
\|h_1\|_1 &= \int_{\bb+\Del^2}^{\infty} \int_{(-\infty, \bb]\times [(n-1)\del, n\del]} \frp_{n\del-s}(y-x) \beta(dx\, ds) \, dy\\ 
&\le \sqrt{2} e^{-\Del^4/4\del} \del \le e^{-\Del^5/\del} \del .
\end{align*}
Here we have used the first condition on $\Del_0$ and the first property of $\del_0$.
Recall that $\|h\|_1 = \del$. Let $q \in [0,1]$ be such that $\|h_2\|_1 = q\del$. Then
$q\ge 1- e^{-\Del^5/\del}$.
We now claim that
\begin{align}\label{eq:q3}
\bb+\Del^2 = \gamma_g^{4\del} + \Del^2 \le \gamma_w^{\Del}.
\end{align}
Assuming the claim, the proof of \eqref{eq:q1} is now completed as follows. 
Write $\calC_{\Del, q\del}w = w-\tilde h$, where $\|\tilde h\|_1 = q\del$.
Note that $\tilde h$ is supported on the right of $\gamma_w^{\Del}$ while 
$h_2$ is supported on the left of $\bb+\Del^2$. This together with the fact that $\|h_2\|_1 = \|\tilde h\|_1 = q\del$ says that
$$\calC_{\Del, q\del}w = w-\tilde h \LE w-h_2 = w-h+h_1 \LE w-h.$$
Also note that, since $q\del \ge \del - \eps$,
$$\calC_{\Del, \del}w \LE \calC_{\Del, q\del}w \mod \eps.$$
Combining the two we have
$$\calC w = \calC_{\Del, \del}w \LE w-h \mod \eps$$
which gives the statement in \eqref{eq:q1}.

We now prove the claim in \eqref{eq:q3}.
We first show that
\begin{equation}\label{eq:aa8}
\gamma_g^{\Del/3} \ge \gamma_g^{4\del} + 2 \Del^2.
\end{equation}
Note that, by our choice of $\del, \Del$,
$4\del < \Del/6$ (see second property of $\del_0$). Thus it suffices to show that, with $\gamma = \gamma_g^{\Del/6}$
$$ g\cc{\gamma,\gamma + 2\Del^2} \le \Del/6.$$
Since 
$v(t, \cdot)$ is nondecreasing, the left side can be bounded by
$$2\Del^2 g(\gamma + 2\Del^2) \le 2\Del^2 g(\gamma_*+1) \le 2\Del^2c_{v,T} (1+ (\gamma_*+1)^m) \le \Del/6.$$
Here, for the first inequality, we have used the first and third conditions on $\Del_0$, for the second we have used \eqref{eq:cvt2} and for the final inequality we have used the second property of $\Del_0$.
This proves the statement in \eqref{eq:aa8}.

We now show that
\begin{equation}\label{eq:q5}
\gamma_{w}^{2\Del/3} = \gamma_{Sg}^{2\Del/3} \ge \theta \doteq \gamma_g^{\Del/3} - \Del^2.
\end{equation}
Write $g= g_1+ g_2$, where $g_2 = \calC_{\Del/3}g$ and $g_1 = g - g_2$.
Note that $\|g_1\|_1 = \Del/3$.
To prove \eqref{eq:q5} it suffices to show that
\begin{equation}\label{eq:q7}
\calS g\cc{-\infty,\theta} \le 2\Del/3.
\end{equation}
Note that
\begin{align*}
\calS g\cc{-\infty,\theta} &= \calS g_1\cc{-\iy,\theta}
+ \calS g_2\cc{-\iy,\theta}\\
&\le \frac{\Del}{3} + \int_{-\infty}^{\theta} \int_{\R} \frp_{\del}(y-x) g_2(x) dx dy\\ 
&= \frac{\Del}{3} + \int_{-\infty}^{\theta} \int_{\theta+\Del^2}^{\infty} \frp_{\del}(y-x) g_2(x) dx dy.
\end{align*}
Write for $(x,y) \in (\theta+\Del^2, \infty)\times (-\infty, \theta)$
$$
(x-y)^2 \ge \frac{1}{2}(x-\theta)^2 + \frac{1}{2}(x-y)^2.$$
Then, with $p^*_{\del}(x, \cdot)$ the density of $x+ \sqrt{2}W_{\del}$,
\begin{multline*}
\int_{-\infty}^{\theta} \int_{\theta+\Del^2}^{\infty} \frp_{\del}(y-x) g_2(x) dx dy \le 
\sqrt{2} \int_{\theta+\Del^2}^{\infty} g_2(x)\exp\{-(x-\theta)^2/4\del\} 
\int_{\R} p^*_{\del}(x,y) dy\, dx\\ 
\le c_{v,T}\sqrt{2}\int_{\theta+\Del^2}^{\infty} (1+ x^m) \exp\{-(x-\theta)^2/4\del\} dx\\ 
\le c_{v,T}\sqrt{2}\int_{\Del^2}^{\infty} (1+ |x+\theta|^m) \exp\{-x^2/4\del\} dx \\ 
\le c_{v,T}\sqrt{2}\exp\{-\Del^4/(8\del)\}\int_{\Del^2}^{\infty}
(1+ |x+\gamma_*+1|^m) \exp\{-x^2/8\del\} dx
 \le \Del/3,
\end{multline*}
where the last inequality is a consequence of the third condition on $\del_0$.

Then we have that
$$
\calS g\cc{-\infty,\theta} \le \frac{\Del}{3} + \Del/3 = 2\Del/3$$
which proves \eqref{eq:q7} and thus also \eqref{eq:q5}.
Finally
$$\gamma_w^{\Del}  \ge \gamma_w^{2\Del/3} \ge \gamma_g^{\Del/3} - \Del^2  \ge \gamma_g^{4\del} + 2\Del^2 - \Del^2 = \gamma_g^{4\del} +\Del^2 $$
which proves the claim in \eqref{eq:q3} and completes the proof of the lemma.
\qed

\begin{proposition}
\label{prop:q3}
Let $v_0 \in \calV_0^\pol(\R)$.
Fix $0<t_0<T$. Let $\Del_0$ and for $\Del \in (0, \Del_0)$, $\delta_0 = \delta_0(\Del)$ be as in Lemma \ref{lem:q5}. Then for $\Del \in (0, \Del_0)$, $\del \in (0, \del_0)$, $n\in \N$, $n\del \le T$, one has
$$ v_{n\delta}^{(\delta, +)} \LE v_{n\delta}^{(\hat \delta,-)} \mod \Del.$$
\end{proposition}

\proof
We argue by induction. Assume $v_{(n-1)\delta}^{(\delta, +)} \LE v_{(n-1)\delta}^{(\hat \delta,-)} \mod \Del$.
Then, from Lemma \ref{lem:simpest}(c), 
$$\calS_{\del} v_{(n-1)\delta}^{(\delta, +)} \LE \calS_{\del} v_{(n-1)\delta}^{(\hat \delta,-)} \mod \Del.$$
Applying Lemma \ref{lem:q7} we now have
$$\calC_{\del}\calS_{\del} v_{(n-1)\delta}^{(\delta, +)} \LE 
\calC_{\hat \del}\calS_{\del} v_{(n-1)\delta}^{(\hat \delta,-)} \mod \Del ,$$ i.e.
$v_{n\delta}^{(\delta, +)} \LE v_{n\delta}^{(\hat \delta,-)} \mod \Del.$
The result follows.
\qed

We can now prove the main result.

\noi{\bf Proof of Theorem \ref{thm:uniq}.}
In view of the relation \eqref{03} between $v$ and $\beta$,
it suffices to to show that if $(v^i, \beta^i)$, $i=1,2$, are two solutions then $v^1=v^2$. Arguing via contradiction, suppose  that $T> 0$ and $r \in \R$ are such that
$v^1_T\cc{-\iy,r} < v^2_T\cc{-\iy,r}$.
Let $\del_n = T/n$ for $n \in \N$. Fix $0<t_0<T$.
Then by Lemmas \ref{lem:q9} and \ref{lem:q5}, there is a $\Del_0> 0$ and for each fixed $\Del \in (0, \Del_0)$ a $\del_0 = \Del(\del_0)$ such that for all $\del \in (0, \del_0)$
$$v_{n\delta}^{(\delta, +)}\cc{-\iy,r}
\le v^1_{n\delta}\cc{-\iy,r}$$
and
$$v^2_{n\del}\cc{-\iy,r}
\le v_{n\delta}^{(\hat \delta,-)}\cc{-\iy,r} + \ell_{n, \delh},$$
when $n\del \le T$, where $\ell_{n, \delh}$ is as in \eqref{eq:aa4}.

Now fix a $\Del \in (0, \Del_0)$ and choose $k$ large enough so that $T/k \in (0, \del_0(\Del))$.
Then applying the above inequalities to $\Del$, $\del=\del_k = T/k$ and $n=k$
$$
v_{T}^{(\delta_k, +)}\cc{-\iy,r} \le v^1_{T}\cc{-\iy,r} < v^2_T\cc{-\iy,r}
\le v_{T}^{(\Del, \delta_k,-)}\cc{-\iy,r} + \ell_{k, \Del, \del_k}.$$
Also, from Proposition \ref{prop:q3} with $\delta = \delta_k$ and $n=k$,
$$v_{T}^{(\Del, \delta_k,-)}\cc{-\iy,r} \le v_{T}^{(\delta_k, +)}\cc{-\iy,r} + \Del.$$
Thus, recalling the definition of $\ell_{k, \Del, \del_k}$
$$v_{T}^{(\delta_k, +)}\cc{-\iy,r} \le v^1_{T}\cc{-\iy,r}
< v^2_T\cc{-\iy,r}
\le v_{T}^{(\delta_k, +)}\cc{-\iy,r} + \Del + t_0+\del_k + e^{-\Del^5/\delta_k} T.$$
Thus 
$$
0\le v^2_T\cc{-\iy,r} - v^1_{T}\cc{-\iy,r} \le \Del + t_0+\del_k + e^{-\Del^5/\delta_k} T.$$
Sending $k\to\infty$, then $\Del \to 0$ and then $t_0 \to 0$ we get that
$ v^2_T\cc{-\iy,r} - v^1_{T}\cc{-\iy,r}=0$
which is a contradiction, completing the proof.
\qed

\section{Proof of main results}
\label{sec5}

In this section we prove Theorems \ref{th1} and \ref{th2}. The former follows directly from the results of Sections \ref{sec2}--\ref{sec5}.
The latter borrows several ideas from \cite{dem19}.

\noi{\bf Proof of Theorem \ref{th1}.}
Corollary \ref{cor1}, which shows the tightness of $(\mu^n,\beta^n)$, Lemma \ref{lem2}, which shows that any subsequential limit $(\mu,\beta)$ of $(\mu^n,\beta^n)$ gives a solution $(v,\beta)$ to \eqref{MFBP} via $v(\cdot,t)=\bF^{\mu_t}(\cdot)$, and finally Theorem \ref{thm:uniq}, showing uniqueness of solutions to \eqref{MFBP}, prove both parts of Theorem \ref{th1}.
\qed

Now consider Theorem \ref{th2}. Let $v_0 \in \calV_0^\pol(\R)$ and $\mu_0$ such that $v_0 = \bF^{\mu_0}$. Suppose as in Theorem \ref{th2}, that
 $\mu_0(dx)\ge\la_0\,{\rm leb}_{[0,\iy)}(dx)$ where $\la_0>0$. Without loss of generality, we assume that $\la_0<2$. Let $(v,\beta)$ be the unique solution to \eqref{MFBP}, and denote
\[
\sig_t=\inf\{x:v(x,t)>0\}, \qquad t>0, \qquad \sig_0=0.
\]
Also, as before, $\mu_t$ is defined as $v(\cdot, t) = \bF^{\mu_t}(\cdot)$.

\begin{lemma}\label{lem51}
	For each $t>0$ and $[a,b] \subset (\sigma_t, \infty)$, one has $\mu_t(a,b) \ge \la_0(b-a)$.
\end{lemma}

\proof
Arguing by contradiction we assume that there exist
$t\in(0,\iy)$, $\la_-\in(0,\la_0)$, $a\in(\sig_t,b)$ and $b\in(\sig_t,\iy)$, such that $\mu_t(a,b) < \la_-(b-a)$. Consider $\atlas(n,\pi^n)$ where $\pi^n$ is the probability law of a PPP with intensity measure $n\mu_0$. This model satisfies Assumption \ref{ass1} (see Example \ref{ex1} and Proposition \ref{prop01}), and so part (b) of Theorem \ref{th1} applies. Denote by $X^n_i$, $Y^n_i$, $Z^n_i$ and $\mu^n$ the corresponding processes (the named particle system, ranked particle system, gaps, normalized configuration measure, resp.). For this model, $\mu^n_t \to \mu_t$ in probability, and since $\mu_t(\{b\})=0$,
we have
\[
		\liminf_{n\to \infty} \PP(\mu^n_t(a,b] < \la_-(b-a)) >0.
\]
Note that $Y^n_0(t) = \inf\{x \in \R: \mu^n_t(-\infty, x] >0\}$.
Then from the upper semicontinuity in $d_*$ of the quantile function (see \cite{dem19} below equation (3.34)),
$\limsup_{n} \PP(Y^n_0(t) \ge a) =0$. Thus
\[
		\liminf_{n\to \infty} \PP(\mu^n_t(a,b] < \la_-(b-a), Y^n_0(t)<a) >0.
\]
	As follows from Lemma \ref{lem1}, in which condition \eqref{eq:bb1} was verified for $\mu^n$, given $\eta >0$, we can find $M >0$ such that
\[
		\limsup_{n\to \infty} \PP(\mu^n_t(-\infty,b] \ge M) \le \eta.
\]
Combining the above two displays, we can find $M>0$ such that
\[
		\liminf_{n\to \infty} \PP(\mu^n_t(a,b] < \la_-(b-a), Y^n_0(t)<a, \mu^n_t(-\infty, b] <M) >0.
\]
Thus, with $c_n =\lfloor\la_-(b-a)n\rfloor +1$,
\begin{equation}\label{r20}
		\liminf_{n\to\iy} \sum_{l=1}^{\lfloor n M\rfloor +1} \PP \left( \sum_{i=l}^{l+c_n} Z^n_i(t) > (b-a)\right) >0.
\end{equation}

We now apply monotonicity arguments from \cite{AS} and \cite{dem19} to the gaps $Z^n_i(t)$.
Let $\{\hat Z^n_i\}$ denote the gap process in $\atlas(n,\hat\pi^n)$ where $\hat\pi^n$ is the law of a $\ppp_{[0,\iy)}(n\la_0)$. Then $Z^n_i(0)$ are stochastically dominated by $\hat Z^n_i(0)$. Hence by \cite[Corollary 3.10]{AS}, $Z^n_i(t)$ are stochastically dominated $\hat Z^n_i(t)$.
Since $\la_0<2$, it follows from \cite[Proposition 1.7]{dem19} that $\hat Z^n_i(t)$ are stochastically dominated by $\hat Z^n_i(0)$. Combining the two statements, $Z^n_i(t)$ are stochastically dominated by i.i.d.\ $\Exp(n\la_0)$. Then the right-hand side of \eqref{r20} increases upon replacing $\{Z^n_i(t)\}$ with i.i.d.\ $\Exp(n\la_0)$. Therefore, with $N_n \sim \poi(\la_0(b-a)n)$,
	$$\liminf_{n\to\iy} n \PP(N_n \le \la_-(b-a)n)>0.$$
	Since $\la_-<\la_0$, the above conclusion is clearly false. This gives a contradiction and completes the proof of the lemma.
\qed

\begin{lemma}\label{lem52}
$\sigma\in C([0,\iy),\R)$.
\end{lemma}

\proof First we show the continuity on $(0,\iy)$.
Define for $q\ge 0$, $\sigma_t(q) = \inf\{x: v(x,t) >q\}$ (and note that $\sig_t(0)=\sig_t$). Fix $t>0$. By the continuity of $v$ in both variables, for any $q>0$
\begin{equation}\label{eq:113}
\lim_{s \to t} v(\sigma_t(q),s) = v(\sigma_t(q),t) = q.
\end{equation}
Thus there is a $\del = \del(q)>0$ such that whenever $|t-s| \le \del$, $v(\sigma_t(q),s) > q/2$ and consequently, for such $s$, $\sigma_t(q)> \sigma_s$.
It then follows that, for such $s$,
$$|v(\sigma_s(q),s) - v(\sigma_t(q),s)| = \mu_s(\sigma_s(q)\wedge \sigma_t(q), \sigma_s(q)\vee \sigma_t(q) ) \ge \la_0 |\sigma_t(q)-\sigma_s(q)|,$$
where we used the fact that $\mu_s(x,y)\ge \la_0(y-x)$ whenever $\sig_s<x<y$ proved in Lemma \ref{lem51}.
Hence
$$|\sigma_t(q)-\sigma_s(q)| \le \la_0^{-1} |v(\sigma_s(q),s) - v(\sigma_t(q),s)| = \la_0^{-1}|q-v(\sig_t(q),s)|.$$
In view of \eqref{eq:113} we have that $\sigma_s(q) \to \sigma_t(q)$ as $s\to t$ and so $t \mapsto \sig_t(q)$ is continuous for any $q>0$.

Next let $0\le q'\le q$. Then
$$ q-q' = \mu_t(\sigma_t(q'),  \sigma_t(q) ) \ge \la_0|\sigma_t(q) - \sigma_t(q')|,$$
showing that $q \mapsto \sigma_t(q)$ is Lipschitz($\la_0^{-1}$) for any $t>0$. 
Combining this with the continuity of $t\mapsto \sigma_t(q)$ for $q>0$, we now get the continuity of the map $t\mapsto \sigma_t$, $t\in(0,\iy)$.

Next we show that $\sig(0+)=0$, establishing the continuity of $\sig$ on all of $[0,\iy)$.
We first argue that $\limsup_{t\to 0}\sigma(t) \le 0$. Arguing by contradiction, suppose that for some $\al>0$, $\limsup_{t\to 0}\sigma(t) >\al$.
	Then, there exists a sequence $t_k \downarrow 0$ such that for all $k$,
	$\sigma(t_k) >\al$.
	This says that, for each $k$,  $v(\al, t_k)=0$. Then
	\begin{align}\label{eq:344n1}
		0 = v(\al/2, t_k) = \int \frp_{t_k}(y-\al/2) v_0(y)dy - \int_0^{t_k} \frp_{t_k-s}(\sigma_s - \al/2) ds.
	\end{align}
	Note that
	$$\int_0^{t_k} \frp_{t_k-s}(\sigma_s - \al/2) ds \le \left(\frac{2}{\pi}\right)^{1/2} t_k^{1/2}.$$
	
	Also, since $v_0(y)=0$ for $y<0$ and for $y>0$, $v_0(y) \ge v_0(y)-v_0(0) \ge \la_0y$ by assumption,
	\begin{align*}
		 \int \frp_{t_k}(y-\al/2) v_0(y)dy &= \int_{0}^{\infty} \frp_{t_k}(y -\al/2) v_0(y) dy \ge \la_0 \int_{0}^{\infty} \frp_{t_k}(y-\al/2) y dy\\
		 &= \la_0 \E((W_{t_k}+\al/2)^+) \ge \la_0(\al/2 - \E|W_{t_k}|)  \ge \la_0(\al/2 - t_k^{1/2}). 
	\end{align*}
	In the above display $W$ is a standard Brownian motion.
	Using the last two displays in \eqref{eq:344n1} we have, for all $n$,
	\begin{align*}
		0 &= v(\al/2, t_k) = \int \frp_{t_k}(y-\al/2) v_0(y)dy - \int_0^{t_k} \frp_{t_k-s}(\sigma_s - \al/2) ds\\
		&\ge \la_0(\al/2 - t_k^{1/2}) - \left(\frac{2}{\pi}\right)^{1/2} t_k^{1/2}.
	\end{align*}
	However the term on the right side is strictly positive for large $k$ which gives a contradiction and completes the proof of $\limsup_{t\to 0}\sigma(t) \le 0$.
	
	Now we show the complementary inequality $\liminf_{t\to 0}\sigma(t) \ge 0$. Again, arguing by contradiction, suppose for some $c>0$, $\liminf_{t\to 0}\sigma(t)  < -c$. Then there exists  a sequence $t_k \downarrow 0$ such that for all $k$, $\sigma(t_k) < -c$ for all $k$.  From Lemma \ref{lem51},
	$$
	v(\sigma(t_k) + c/4, t_k) =  v(\sigma(t_k) + c/4, t_k)  - v(\sigma(t_k), t_k)  \ge \la_0 c/4.$$
	Also, recalling that $v_0(y) = 0$ for $y<0$ and for some $c_1>0$, $v_0(y) \le c_1(1+y^{c_1})$ for $y\ge 0$,
	\begin{align*}
		\la_0 c/4 &\le v(\sigma(t_k) + c/4, t_k) \le \int \frp_{t_k}(y-(\sig(t_k) +c/4)) v_0(y)dy \\
		&= \int_0^{\infty} \frp_{t_k}(y-(\sig(t_k) +c/4)) v_0(y)dy\\
		&\le \int_0^{\infty} \frac{c_1}{(2\pi t_k)^{1/2}}  e^{-(y+c/2)^2/2t_n} (1+y^{c_1})dy\\
		&\le c_1e^{-c^2/16t_k} \int_{\R} \frac{1}{(2\pi t_k)^{1/2}}  e^{-(y+c/2)^2/4t_k} (1+|y|^{c_1})dy \le c_2 e^{-c^2/16t_k}
		\end{align*}
		for some $c_2>0$. Sending $k\to \infty$ we arrive at a contradiction. This completes the proof of $\liminf_{t\to 0}\sigma(t) \ge 0$.
\qed

\begin{lemma}\label{lem53}
	We have $\beta(dx, dt ) = \delta_{\sig_t}(dx) dt$.
\end{lemma}

\proof
Denote
\[
G=\{(x,t):x<\sig_t,t>0\},
\qquad
H=\{(x,t):x>\sig_t,t>0\}.
\]
By Lemma \ref{lem52}, both are open subsets of $\R^2$. In view of the representation $\beta(dx,dt)=\beta_t(dx)dt$,
it suffices to prove that $\beta(G)=\beta(H)=0$.

To show that $\beta(H)=0$, note by monotonicity of $x\mapsto v(x,t)$ that $v>0$ on $H$. Hence the statement follows from $\int v\,d\beta=0$.

It remains to show that $\beta(G)=0$. Arguing by contradiction, assume $\beta(G)>0$. Then there exists a rectangle $R=[x_0,y_0]\times[s_0,t_0]\subset G$ and $z>y_0$ such that $\beta(R)>0$ and $\sig_s>z$ for all $s\in[s_0,t_0]$. We show that this contradicts nonnegativity of $v$.

Denote by $m$ the measure on $[s_0,t_0]$ defined via $m[s,t]=\beta([x_0,y_0]\times[s,t])$, for $s_0 \le s \le t \le t_0$. By assumption, $S:=m[s_0,t_0]>0$. Denoting $\eps=S(t_0-s_0)^{-1}/2$, it follows by a pigeonhole argument that for every $\del\in(0,t_0-s_0)$ there is $s\in[s_0,t_0-\del]$ with $m[s,s+\del]\ge\eps\del$. Now, by \eqref{83}, for $s_0\le s<t\le t_0$, $\del=t-s$,
\[
v(x,t)=\calS_\del v_s(x)-\int_s^t\calS_{t-\theta}\beta_\theta(x)d\theta,
\qquad x\in\R,
\]
where $v_s=v(\cdot,s)$. Also,
\begin{align*}
\int_{-\iy}^{y_0}\int_{s}^{t}\calS_{t-\theta}\beta_\theta(x)d\theta dx
&\ge\int_{-\iy}^{y_0}\int_{s}^{t}\int_{[x_0,y_0]}\frp_{t-\theta}(y-x)\beta_\theta(dy)d\theta dx
\\
&\ge\frac{1}{2}\int_{s}^{t}\int_{[x_0,y_0]}\beta_\theta(dy)d\theta
=\frac{1}{2}m[s,t] \ge \frac{1}{2}\eps\del.
\end{align*}
Moreover,  for $c$ not depending on $\del$,
\begin{align*}
\int_{-\iy}^{y_0}\calS_{\del}v_{s}(x)dx
&=\int_{-\iy}^{y_0}\int_\R \frp_{\del}(y-x) v_{s}(y)dydx
\\
&\le\int_{-\iy}^{y_0}\int_z^\iy \frp_{\del}(y-x) c(1+y_+^c)dydx.
\end{align*}
With $a=z-y_0$, note that in the above integrand, $y-x\ge a$. Hence the right-hand side above is bounded by
$c_1e^{-c_2a^2/\del}$. As a result,
\[
\int_{-\iy}^{y_0}v(x,t)dx\le c_1e^{-c_2 a^2/\del}-\frac{1}{2}\eps\del.
\]
For $\del>0$ small, the above is negative, a contradiction. This completes the proof of the lemma.
\qed

\noi{\bf Proof of Theorem \ref{th2}.}
Part (a) follows from Lemmas \ref{lem52} and \ref{lem53}. Part (b) is immediate from part (a) and Definitions \ref{def00} and \ref{def0}.

Consider now part (c). Thus we suppose that Assumption \ref{ass1} holds and $v_0 = \bF^{\mu_0}$.
Fix $T>0$ and let
$$\gamma_n \doteq \sup_{t\in[0,T]}(Y^n_0(t)-\sig_t)^+,
\qquad
\del_n=\sup_{t\in[0,T]}(Y^n_0(t)-\sig_t)^-.
$$
It suffices to show that $\gamma_n, \delta_n$ converge to $0$ in probability as $n\to \infty$.
Consider first $\gamma_n$.
We invoke Skorohod's representation, by which we may assume that there exists a full-measure event $\Om_1$ such that for all $\om\in\Om_1$, $\mu^n\to\mu$ in $C(\R_+,\calM_*)$. It suffices to show that for all $\om \in \Om_1$, $\limsup_{n\to \infty} \gamma_n(\om) \doteq \gamma =0$. Arguing by contradiction, suppose that $\gamma>0$. Then (suppressing $\om$), there exists a sequence $(t_k, n_k)_{k \in \N}$ such that $t_k \in [0,T]$ for all $k$, $n_k \uparrow \infty$ as $k \to \infty$, and, for all $k$,
$$Y^{n_k}_0(t_k) \ge \sigma(t_k) + \gamma/2.$$
We can assume without loss of generality that $t_k \to t$ for some $t \in [0,T]$.
Since $\mu^{n_k}(t_k) \to \mu(t)$ as $k\to \infty$, the upper semicontinuity of the quantile function, that was mentioned in the proof of Lemma \ref{lem51} and the property that $\sigma(t_k) \to \sigma(t)$, implies that 
$$\sigma(t) + \gamma/2 = \lim_{k\to \infty} \sigma(t_k) + \gamma/2 \le   \limsup_{k\to\iy}Y^{n_k}_0(t_k)\le\sig(t).$$
This gives a contradiction and proves that $\gamma =0$.

We now argue that $\delta_n \to 0$ in probability. Recall the constant $c$ from the statement of the theorem. Denote by $\tilde \pi_n$ the law of a PPP on $\R_+$ with intensity measure $ncdx$.
Using strong approximative solutions of the infinite Atlas model (see \cite[Corollary 3.10]{AS}), we can construct gap sequence processes $Z^n, \tilde Z^n$ of $\atlas(n,\pi_n)$ and $\atlas(n,\tilde \pi_n)$ on a common probability space, with $Y^n_0(0) = \tilde Y^n_0(0) =0$ (i.e. the lowest particle starting at $0$ in both models), such that
\begin{equation}\label{eq:ab228}
Z^n(t) \le \tilde Z^n(t) \mbox{ for all } t \ge 0.\end{equation}
The normalized configuration measure processes for the two systems will be denoted by $\mu^n$ and $\tilde \mu^n$.
For rest of the proof we work on this probability space.  By \cite[Theorem 1.2]{dem19}, $(\tilde Y^n_0(t), \tilde\mu^n_t) \to (\tilde\sig_t, \tilde\mu_t)$, uniformly on compacts, in probability, as $n\to \infty$. Furthermore, $\tilde\mu_t$ has density bounded away from $0$ to the right of $\tilde\sig_t$, in particular $$\inf_{t \in [0,T]}\tilde\mu_t[\tilde\sig_t,\tilde\sig_t+\al]>a\al, \; \al \in (0,1],$$
for some $a>0$.
These two facts imply that
\begin{equation}\label{eq:ab206}\limsup_{n\to \infty} \PP(\inf_{t\in [0,T]} \tilde \mu^n_t[\tilde Y^n_0(t), \tilde Y^n_0(t)+\al] \le a\al/2) =0, \; \al \in (0,1].
\end{equation}

Arguing by contradiction again, suppose for some $\del>0$,
$$\limsup_{n\to \infty}\PP(\delta_n >\delta) \doteq \theta >0.
$$
We now argue that
\begin{equation}\label{eq:cl1}
\sup_{t \in [0,T]} \mu^n_t(-\infty, \sigma_t - \del/8] \to 0 \mbox{ in probability, as } n \to \infty.
\end{equation}
To see this, let $\bar \sigma \doteq \sup_{t \in [0,T]} \sig_t$, and for each $t \in [0,T]$, define 
continuous, non-increasing function $f^t$ with support in $(-\infty, \bar \sig]$  such that $f^t = 1$ on $(-\infty, \sigma_t - \del/8]$, and $0$ on $[\sig_t, \infty)$, and such that the Lipschitz constant is  $8\delta^{-1}$.
From convergence of $\mu^n$ to $\mu$ we have that
$$\sup_{t \in [0,T]} |\lan \mu^n_t, f^t\ran - \lan \mu_t, f^t\ran | \to 0 \mbox{ in probability, as } n \to \infty.$$
The claim in \eqref{eq:cl1} now follows on noting that, for all $t \in [0,T]$, $\lan \mu^n_t, f^t\ran  \ge \mu^n_t(-\infty, \sigma_t - \del/8]$ and $\lan \mu_t, f^t\ran=0$.
We now have that, with $M^n_t = n \mu^n_t(-\infty, \sigma_t - \del/8]$, and for arbitrary $\eps>0$,
\begin{equation}\label{eq:ab229}
\limsup_{n\to \infty} \PP\left(\sup_{t \in [0,T]} M^n_t \le n \eps, \, \delta_n > \delta\right) = \theta.
\end{equation}
For $\om \in A_n \doteq \{\sup_{t \in [0,T]} M^n_t \le n \eps, \, \delta_n > \delta\}$ (and suppressing $\om$ from the notation),
$Y^n_0(t) \le \sig_t - \del$ for some $t \in [0,T]$. Recalling the definition of $M^n_t$, for such $t$,
$$\del/2 \le \sum_{i=1}^{\lfloor M^n_t\rfloor+1} Z^n_i(t) \le \sum_{i=1}^{\lfloor n\eps\rfloor+1} Z^n_i(t)
\le \sum_{i=1}^{\lfloor n\eps\rfloor+1} \tilde Z^n_i(t),$$
where the last inequality uses \eqref{eq:ab228}.
This says that, on $A^n$,
$$\inf_{t\in [0,T]} \tilde \mu^n_t[\tilde Y^n_0(t), \tilde Y^n_0(t)+\del/2] \le \eps + n^{-1}$$
and consequently
$$
\limsup_{n\to \infty} \PP(\inf_{t\in [0,T]} \tilde \mu^n_t[\tilde Y^n_0(t), \tilde Y^n_0(t)+\del/2] \le 2\eps) \ge \limsup_{n\to \infty} \PP(\inf_{t\in [0,T]} \tilde \mu^n_t[\tilde Y^n_0(t), \tilde Y^n_0(t)+\del/2] \le \eps + n^{-1}) \ge \theta,
$$
where the last inequality uses \eqref{eq:ab229}.
Taking $\eps = a/4$ in the last display and applying \eqref{eq:ab206} with $\al=1$ we arrive at a contradiction, completing the proof of $\delta_n \to 0$ in probability.
\qed

\appendix

\section{Classical solutions to \eqref{FBP}}
\label{sec:app1}

Although we do not deal with classical solutions to \eqref{FBP} in this paper, we state here, for completeness, their relation to weak solutions of \eqref{FBP} and to solutions of \eqref{MFBP} (see e.g.\ \cite{fried68} for a well known similar statement regarding classical and weak solutions to the Stefan problem).
A classical solution to \eqref{FBP} is a pair $(\sig,v)$,
satisfying $\sig\in C((0,\iy),\R)$, $\sig_{0+}=0$, and
$v\in C^{2,1}(\{(x,t):x>\sig_t,t>0\})$, where $v$ and $\pl_x v$ have continuous extensions to $\{(x,t):x\ge\sig_t\}$, and \eqref{FBP} holds.
\begin{proposition}\label{prop0}
Assume that the initial data $v_0$ is continuous, vanishes in $(-\iy,0]$, positive in $(0,\iy)$,
let $(v,\sig)$ be a classical solution
of \eqref{FBP}, and assume that $\pl_tv$ extends continuously to $\{(x,t):x\ge\sig_t,t\ge0\}$. Extend $v$ to $\R\times\R_+$ by $v(x,t)=0$ for $x<\sig_t$, and assume that $v \in C^\pol_\loc(\R\times(0,\iy))$. Then $(v,\sig)$ is a weak solution to \eqref{FBP} in the sense of Definition \ref{def00}, and gives a solution $(v,\beta)$ to \eqref{MFBP} in the sense of Definition \ref{def0} via $\beta(dx,dt)=\del_{\sig_t}(dx)dt$.
\end{proposition}

\vspace{-.5em}
\proof
This follows by a direct calculation.
\qed

\subsection*{Acknowledgments}
RA was partially supported by ISF (grant 1035/20).
AB was partially supported by NSF DMS-2152577, NSF  DMS-2134107, NSF DMS-2506010.

\vspace{-.5em}

\footnotesize

\bibliographystyle{is-abbrv}

\bibliography{main}

%
%

\end{document}